\documentclass[11pt]{amsart}

\usepackage{amsaddr}
\usepackage{amsmath,amssymb,amscd,amsthm}

\usepackage[dvips]{graphicx}
\usepackage{epsfig}
\usepackage{graphics}
\usepackage{color}

\usepackage{url}

\newtheorem{dfn}{Definition}[section]
\newtheorem{thm}[dfn]{Theorem}

\newtheorem{lem}[dfn]{Lemma}
\newtheorem{cor}[dfn]{Corollary}
\newtheorem{conj}{Conjecture}

\theoremstyle{definition}
\newtheorem{rem}[dfn]{Remark}

\newtheorem{prb}{Problem}

\theoremstyle{plain}

\numberwithin{equation}{section}

\def\R{{\mathbb R}}

\def\Sph{{\mathbb S}}

\def\Euc{{\mathbb E}}

\def\H{{\mathbb H}}

\def\dS{\mathrm{d}{\mathbb S}}

\def\phi{\varphi}
\def\epsilon{\varepsilon}
\def\emptyset{\varnothing}

\def\P{{\mathcal P}}
\def\Q{{\mathcal Q}}

\renewcommand{\span}{\operatorname{span}}
\newcommand{\area}{\operatorname{area}}

\def\Vol{\operatorname{Vol}}
\def\Area{\operatorname{Area}}

\def\HE{\operatorname{HE}}

\newcommand{\bigdot}{{\displaystyle{\cdot}}}

\title[Rigidity of surfaces and the Hilbert-Einstein functional I]{Infinitesimal rigidity of convex surfaces through the second derivative of the Hilbert-Einstein functional\\ I: Polyhedral case}
\author{Ivan Izmestiev}
\address{Institut f\"ur Mathematik, MA 8-3 \\
Technische Universit\"at Berlin \\
Stra{\ss}e des 17. Juni 136 \\
D-10623 Berlin, Germany}
\email{izmestiev@math.tu-berlin.de}
\thanks{Supported by the DFG Research Unit 565 ``Polyhedral Surfaces''}
\date{May 25, 2011}

\begin{document}

\begin{abstract}
The paper is centered around a new proof of the infinitesimal rigidity of convex polyhedra. The proof is based on studying derivatives of the discrete Hilbert-Einstein functional on the space of ``warped polyhedra'' with a fixed metric on the boundary. This approach is in a sense dual to using derivatives of the volume in order to prove the Gauss infinitesimal rigidity of convex polyhedra, which deals with deformations that preserve face normals and face areas.

In the spherical and in the hyperbolic-de Sitter space, there is a perfect duality between the Hilbert-Einstein functional and the volume, as well as between both kinds of rigidity.

We also discuss directions for future research, including elementary proofs of the infinitesimal rigidity of hyperbolic (cone-)manifolds and development of a discrete Bochner technique.
\end{abstract}

\maketitle

\setcounter{tocdepth}{1}
\tableofcontents

\section{Introduction}
\subsection{Infinitesimal rigidity of convex polyhedra}

Among various concepts of rigidity, the one that we deal with here is the \emph{infinitesimal rigidity}. Let $P \subset \R^3$ be a compact convex polyhedron; assume for simplicity that all faces of $P$ are triangles. Then an \emph{infinitesimal deformation} of $P$ is an assignment of a vector $q_i$ to each vertex $p_i$. Let $P(t)$ be a family of polyhedra with the same combinatorics as $P$ and vertices $p_i(t) = p_i + t q_i$. An infinitesimal deformation is called \emph{isometric}, if the edge lengths of $P(t)$ remain constant in the first order of $t$ at $t=0$. There always exist \emph{trivial} isometric infinitesimal deformations that are restrictions of infinitesimal isometries of~$\R^3$.
\begin{quote}
Every convex polyhedron $P$ is infinitesimally rigid, that is every isometric infinitesimal deformation of $P$ is trivial.
\end{quote}
This theorem was first stated by Dehn in \cite{Dehn16}. Legendre-Cauchy's argument from \cite{Cau13} easily carries over to the infinitesimal setting. Other proofs are given in \cite{Ale05, Fil92, Pak06, Scl07, Tru81, Weyl17}.

Yet another proof is presented in this paper. In some aspects it is similar to Schlenker's proof \cite{Scl07} and to the Pogorelov-Volkov proof of global rigidity \cite{Pog56, Vol55}. The actual goal of this paper is to describe the general framework into which this argument fits and to indicate possible further developments.

\subsection{The approach}
\label{subsec:IntroAppr}
Instead of deforming the embedding $P \subset \R^3$ by moving the vertices, we deform the metric in the interior of $P$. We choose a point $p_0 \in \mathrm{int}\, P$ and subdivide $P$ into pyramids with $p_0$ as the apex and faces of $P$ as bases. By denoting $r_i = \|p_0 - p_i\|$, we start to vary lengths $r_i$, while leaving the lengths of boundary edges constant. As a result, the total dihedral angles $\omega_i$ around interior edges $p_0p_i$ may become different from $2\pi$, so that the polyhedron $P$ becomes what we call a \emph{warped polyhedron}.

Can one vary the lengths $r_i$ so that the angles $\omega_i$ remain constant in the first order? One can always do this by moving the point $r_0$ inside $P$. It turns out that if these are the only possibilities, then the polyhedron $P$ is infinitesimally rigid: infinitesimal rigidity of $P$ is equivalent to
\begin{equation}
\label{eqn:IntroDKDR}
\dim \ker \left( \frac{\partial \kappa_i}{\partial r_j} \right) = 3,
\end{equation}
where $\kappa_i = 2\pi - \omega_i$ is the \emph{curvature} of the edge $p_0p_i$. Our proof of the infinitesimal rigidity of convex polyhedra goes by determining the rank of the above Jacobian.

An immediate generalization of this approach is to consider an arbitrary triangulation of the polyhedron $P$ and the Jacobi matrix $\left( \frac{\partial \kappa_{ij}}{\partial r_{kl}} \right)$, where $r_{ij}$ and $\kappa_{ij}$ are the length, respectively, curvature of the edge $p_i p_j$, and among the points $p_i$ some are vertices of $P$, some lie on its edges or faces, and some in the interior. (As before, lengths of boundary edges are assumed constant.) Then infinitesimal rigidity of $P$ is equivalent to
$$
\dim \ker \left( \frac{\partial \kappa_{ij}}{\partial r_{kl}} \right) = 3m + n,
$$
where $m$ is the number of interior vertices, and $n$ is the number of vertices interior to faces of $P$. In our joint paper with Jean-Marc Schlenker \cite{IS10}, we have shown that the above Jacobian has exactly $m$ positive eigenvalues. In particular, the matrix has negative spectrum if all points $p_i$ lie at vertices or on the edges of $P$.

\subsection{Hyperbolic manifolds with boundary and hyperbolic cone-manifolds}
\label{subsec:HypMan}
Similar to the previous subsection, consider a triangulation of a compact closed hyperbolic 3--manifold and the Jacobi matrix $\left( \frac{\partial \kappa_{ij}}{\partial r_{kl}} \right)$ of edge curvatures with respect to edge lengths.

\begin{conj}
\label{conj:CalWeil}
The matrix $\left( \frac{\partial \kappa_{ij}}{\partial r_{kl}} \right)$ has corank $3n$ and exactly $n$ positive eigenvalues, where $n$ is the number of vertices in the triangulation.
\end{conj}

The first part, corank equals $3n$, is known. It follows from the infinitesimal rigidity of compact closed hyperbolic manifolds in dimension greater than $2$, the so-called Calabi-Weil rigidity, \cite{Cal61, Wei60}. Therefore, if the arguments used in this paper can be extended to prove Conjecture \ref{conj:CalWeil}, this would yield an elementary proof of Calabi-Weil rigidity. Besides, such extension should lead to a proof of the following conjecture, communicated to me by Jean-Marc Schlenker who proved a smooth analog in \cite{Scl06}.

\begin{conj}
Compact hyperbolic manifolds with convex polyhedral boundary are infinitesimally rigid.
\end{conj}

More generally, consider a compact closed manifold glued from hyperbolic simplices so that the total dihedral angles around the edges may be different from $2\pi$. It is called a (triangulated) hyperbolic \emph{cone-manifold}. Due to works of Hodgson-Kerckhoff \cite{HK98}, Mazzeo-Montcouquiol \cite{MM09}, and Hartmut Weiss \cite{Wei09}, it is known that compact closed hyperbolic manifolds with cone angles less than $2\pi$ are infinitesimally rigid, that is cannot be deformed so that their cone angles remain constant in the first order. Again, this can be reformulated in terms of the rank of the Jacobian of the map $r \mapsto \kappa$, and a generalization of our method would yield an alternative, elementary proof.

Note that the condition on cone angles is similar to the convexity condition for polyhedra. There exist infinitesimally flexible non-convex polyhedra as well as infinitesimally flexible cone-manifolds with some cone angles greater than $2\pi$, see \cite{Cas98, Izm11}.

A study of deformations of Euclidean cone-manifolds from a similar point of view is undertaken in \cite{Gli11}.

\subsection{Discrete Hilbert-Einstein functional}
\label{subsec:IntroDHE}
Here we sketch our proof of $\dim \ker \left( \frac{\partial \kappa_i}{\partial r_j} \right) = 3$ for a star-like triangulation of a convex polyhedron $P$. Consider the function
$$
\HE(r) = \sum_i r_i \kappa_i + \sum_{\{i,j\}} \ell_{ij} \lambda_{ij},
$$
where $\ell_{ij}$ is the length of, and $\lambda_{ij}$ is the exterior dihedral angle at an edge $p_ip_j$ of $P$. We call it the \emph{discrete Hilbert-Einstein functional}, as it is the discrete analog of twice the total scalar curvature of $P$ plus half of the total mean curvature of $\partial P$. The Schl\"afli formula implies that $\frac{\partial \HE}{\partial r_i} = \kappa_i$, therefore
$$
\frac{\partial \kappa_i}{\partial r_j} = \frac{\partial^2 \HE}{\partial r_i \partial r_j},
$$
in particular the left hand side is symmetric in $i$ and $j$. This allows to reformulate the infinitesimal rigidity in the following way:
\begin{equation}
\label{eqn:ReformHEDDot}
\mbox{If }\HE^{\bigdot\bigdot} = 0, \mbox{ then }\dot r \mbox{ is a trivial variation of }r.
\end{equation}
Here $\HE^{\bigdot\bigdot} = \sum_{i,j} \frac{\partial^2 \HE}{\partial r_i \partial r_j} \dot r_i \dot r_j = \sum_i \dot \kappa_i \dot r_i$ is the second derivative of $\HE$ in the direction $\dot r$, and a trivial variation can be defined as the one that preserves boundary dihedral angles: $\dot \lambda_{ij} = 0$.

In order to prove \eqref{eqn:ReformHEDDot}, we derive formula \eqref{eqn:SDDot2} for $\HE^{\bigdot\bigdot}$. This formula implies that $\HE^{\bigdot\bigdot}$ is non-positive, and vanishes only if $\dot r$ is trivial. But as $\HE^{\bigdot\bigdot}$ vanishes by assumption, the variation $\dot r$ must be trivial.

This argument is reminiscent of Koiso's proof \cite{Koi78} of the infinitesimal rigidity of Einstein manifolds under certain assumptions on the curvature operator. Koiso uses integration by parts to obtain two formulas for the second derivative of the Hilbert-Einstein functional. This yields an equation with zero on one side, while on the other side one has a non-positive quantity that vanishes only if the deformation is trivial. Note that our formula \eqref{eqn:SDDot2} is obtained from $\HE^{\bigdot\bigdot} = \sum_i \dot \kappa_i \dot r_i$ by kind of discrete integration by parts.

Koiso's proof is an example of application of Bochner's technique: a second order differential operator is expressed as the sum of Laplacian and of a non-negative 0--th order operator. Thus our proof should be a particular manifestation of a \emph{discrete Bochner technique}.

\subsection{Discrete Bochner technique}
\label{subsec:DiscrBochner}
Einstein manifolds in dimension $3$ are manifolds of constant sectional curvature, therefore Koiso's theorem contains Calabi-Weil rigidity of hyperbolic $3$--manifolds as a special case. Weil also uses Bochner's technique; his approach is related to that of Koiso in the same way as moving vertices of a polyhedron $P$ is related to deforming the metric inside $P$. The infinitesimal rigidity of hyperbolic cone-manifolds with cone angles less than $2\pi$ (see Subsection \ref{subsec:HypMan}) is proved by extending Weil's arguments with the help of Cheeger's Hodge theory for singular spaces.

It will be just natural if the infinitesimal rigidity of cone-manifolds can be reproved by developing a discrete Bochner technique for this situation. Such a proof would not only be elementary, it would also provide a discrete-geometric counterpart to the original argument of Hodgson, Kerckhoff, and others.

Another theorem where the need for a discrete Bochner technique is felt, is Cheeger's discrete analog \cite{Che86} of Bochner-Gallot-Meyer vanishing theorem:
\begin{quote}
If in a Euclidean cone-manifold of dimension $d$ all cone angles are less than $2\pi$, then it is a real homology $d$--sphere.
\end{quote}
Again, the known proof uses Hodge theory for singular spaces.

It should be mentioned that Forman developed a \emph{combinatorial Bochner technique} \cite{For03} that takes into account the combinatorics of a simplicial complex (or, more generally, cell complex), but not its geometry, expressed by edge lengths.

One of the components for a discrete Bochner technique should be the discrete Hodge theory based on Whitney forms, see e.~g. \cite{Dod76}. The other component presumes some sort of discrete Riemannian geometry, still to be found.

\subsection{Volume derivatives and the Alexandrov-Fenchel inequality}
\label{subsec:IntroVolDer}
If one goes into details of the proof of $\dim \ker \left( \frac{\partial \kappa_i}{\partial r_j} \right) = 3$ sketched in Subsection \ref{subsec:IntroDHE}, then one sees that it resembles very much a known proof of
\begin{equation}
\label{eqn:RankDA}
\dim \ker \left( \frac{\partial A_i}{\partial h_j} \right) = 3,
\end{equation}
where $h_i$ are lengths of perpendiculars dropped to faces of a convex polyhedron $Q$ from the origin, $A_i$ are face areas, and the polyhedron is deformed by varying $h_i$ while keeping the directions of face normals fixed.

In fact, equation \eqref{eqn:RankDA}, with $d$ on the right hand side, also holds for convex $d$--dimensional polyhedra. This is a key lemma in the proof of the Alexandrov-Fenchel inequalities. Alexandrov-Fenchel inequalities also imply that the Jacobian in \eqref{eqn:RankDA} has exactly one positive eigenvalue (cf. the end of Subsection \ref{subsec:IntroAppr}) and describe a part of the positive cone.

It is geometrically clear that $A_i = \frac{\partial \Vol}{\partial h_i}$, where $\Vol(h)$ is the volume of a polyhedron with given face normals and with support numbers $h_i$. Therefore equation \eqref{eqn:RankDA} computes the rank of the second derivative of the volume, thus completing the analogy with Subsection \ref{subsec:IntroDHE}.

\subsection{Duality between the volume and Hilbert-Einstein functional}
The analogy described in the previous subsection culminates in a striking identity
$$
\frac{\partial^2 \HE}{\partial r_i \partial r_j} = \frac{\partial^2 \Vol}{\partial h_i \partial h_j},
$$
where the polyhedron $Q$ on the right hand side is polar dual to the polyhedron $P$ on the left hand side. We were able to prove the above identity only by a direct computation. Note that functions $\HE(r)$ and $\Vol(h)$ have different nature: the former uses inverse trigonometric functions while the latter is polynomial.

The situation is nicer in spherical and hyperbolic geometry. For a spherical convex polyhedron $P$ and its polar dual $P^*$ we have by \cite{McM75}
$$
\Vol(P) + \frac{1}{2} \sum_{\{i,j\}} \ell_{ij} \lambda_{ij} + \Vol(P^*) = \pi^2.
$$
This implies $S(P) + S^*(P^*) = 2\pi^2$, where $S$ and $S^*$ are functionals on the space of warped spherical polyhedra. Functionals $S$ and $S^*$ are the true analogs of $\HE$ and $\Vol$, respectively; in particular their variational properties are similar to those of their Euclidean counterparts. For more details, see Subsections \ref{subsec:SpherDual} and \ref{subsec:SpherDualMetr}. The hyperbolic case is similar, but the polar dual $P^*$ of a hyperbolic polyhedron $P$ lives in the de Sitter space, see Subsection \ref{subsec:HypDS}.

\subsection{Gauss infinitesimal rigidity and Minkowski theorem}
While equation \eqref{eqn:IntroDKDR} is equivalent to the infinitesimal rigidity of a polyhedron, equation \eqref{eqn:RankDA} implies what we call \emph{Gauss infinitesimal rigidity}:
\begin{quote}
If the support numbers $h_i$ of a convex polyhedron vary in such a way that the face areas remain constant in the first order, then the polyhedron undergoes a parallel translation.
\end{quote}
This theorem was stated and proved by Alexandrov \cite[Chapter XI]{Ale05}. Minkow\-ski proved what we would call global Gauss rigidity: two convex polyhedra with the same face normals and face areas differ by a parallel translations. Existence is asserted in the \emph{Minkowski theorem}:
\begin{quote}
Given unit vectors $\nu_i$ that span $\R^3$ and positive numbers $C_i$ such that $\sum_i C_i \nu_i = 0$, there exists a convex polyhedron with face areas $C_i$ and outward face normals $\nu_i$.
\end{quote}
More generally, Minkowski existence and uniqueness theorem holds in $\R^d$.

\subsection{Alexandrov theorem on existence of a convex polyhedron with a given metric on the boundary}
A counterpart to the Minkowski theorem in dimension $3$ is the \emph{Alexandrov theorem} \cite{Ale42}:
\begin{quote}
Given a Euclidean cone-metric $g$ on the 2--sphere $\Sph^2$ with all cone angles less than $2\pi$, there exists a convex polyhedron in $\R^3$ with $g$ as the intrinsic metric on the boundary.
\end{quote}
In a joint paper \cite{BI08} with Alexander Bobenko, we gave a new proof of the Alexandrov theorem. The polyhedron $P$ is obtained by constructing a family of warped polyhedra $P_t, 0 \le t \le 1$ where the curvatures $\kappa_i$ tend to $0$ as $t$ tends to $1$. The local existence of such family is based on the following property of the Jacobian of the map $r \mapsto \kappa$:
$$
\dim \ker \left( \frac{\partial \kappa_i}{\partial r_j} \right) = 0, \mbox{ if } 0 < \kappa_i < \delta_i \mbox{ for all } i.
$$
Here $\delta_i$ is the angular defect of the $i$-th cone point on $(\Sph^2, g)$. For the global existence one has to make sure that polyhedra don't degenerate in the process of deformation. Note that the triangulation of the boundary of $P$ may change in the process of deformation. See also Subsection \ref{subsec:ExistThms}.

\subsection{Acknowledgements}
I would like to thank Fran\c{c}ois Fillastre and Jean-Marc Schlenker for interesting discussions and useful remarks.

\section{Gauss rigidity of convex polyhedra}
\label{sec:GaussRig}
\subsection{The theorem}
\label{subsec:Theorem}
A set $(Q_t)_{0 \le t < \epsilon}$ of compact convex polyhedra is called a \emph{linear family}, if each $Q_t$ is obtained from $Q = Q_0$ by parallelly translating the planes of the faces, with the translation vector of each face depending linearly on $t$.

In this section, a proof of the following theorem is presented.
\begin{thm}
\label{thm:MinkRigDiscr}
Assume that the area of each face in the linear family $(Q_t)$ is constant in the first order of
$t$. Then all $Q_t$ are translates of $Q$.
\end{thm}

This theorem is proved by Alexandrov in \cite[Chapter XI]{Ale05} by two different methods. The proof given here is essentially the second proof of Alexandrov, but presented in a self-consistent way. Our purpose is to reveal that the proof is based on certain variational properties of the volume, in order to make the relationship with Section \ref{sec:MetrRig} more straightforward.

Let us introduce some notations. Let $F_1, F_2, \ldots, F_n$ be the faces of $Q$, and let $\nu_i$ be the outer unit normal to the face $F_i$. Then the equation of the plane of the face $F_i$ is
$$
\span(F_i) = \{x \in \R^3\ |\ \langle x, \nu_i \rangle = h^0_i\},
$$
where $h^0_i$ is the signed distance from the coordinate origin $0 \in \R^3$ to the plane spanned by $F_i$. The numbers $(h_i^0)$ are called \emph{support parameters} of the polyhedron $Q$. We have
$$
Q = \{x \in \R^3\ |\ \langle x, \nu_i \rangle \le h_i^0,\ i=1,\ldots,n\} =: Q(h^0).
$$

We fix the directions of outer normals, and let the support parameters vary. A \emph{linear family}
of polyhedra is $Q(h^0+tu)$. A polyhedron $Q(h)$ with $h_i = h_i^0 + \langle a, \nu_i
\rangle$ is a parallel translate of $Q$ by $a \in \R^3$.

Let $A_i$ be the area of the face $F_i$. Denote by $DA_i(u)$ the derivative of $A_i$ in the
direction $u$:
$$
D A_i (u) := \left.\frac{d}{dt}\right|_{t=0} A_i(h+tu) = \sum_i \frac{\partial A_i}{\partial h_i} u_i.
$$
We will usually put $\dot h$ in place of $u$ and then write
$$
D A_i (\dot h) =: \dot A_i.
$$
Thus the statement of Theorem \ref{thm:MinkRigDiscr} can be rewritten as

\begin{quote}
\emph{If $\dot h \in \R^n$ is such that at $h=h^0$ we have $\dot A_i = 0$ for all $i$, then there
exists $a\in \R^3$ such that $\dot h_i = \langle a, \nu_i \rangle$ for all $i$.}
\end{quote}

It is easy to see that the space of ``trivial deformations'' $\dot h_i = \langle a, \nu_i \rangle$
has dimension $3$. Hence the theorem can also be reformulated as

\begin{quote}
\emph{The Jacobian of the map $(h_i) \mapsto (A_i)$ has corank $3$:
$$
\dim \ker \left( \frac{\partial A_i}{\partial h_j} \right)\bigg|_{h=h^0} = 3.
$$
}
\end{quote}

\subsection{The approach}
\label{subsec:Approach}
Recall that we put
$$
Q(h) := \{x \in \R^3\ |\ \langle x, \nu_i \rangle \le h_i,\ i=1,\ldots,n\}.
$$
The set $Q(h)$ is a convex polyhedron for all $h \in \R^n$, but it may have less than $n$ faces or
even be empty. Let $U \subset \R^n$ be a neighborhood of $h^0$ such that $Q(h)$ has $n$ faces for
all $h \in U$. Consider the function
$$
\Vol \colon U \to \R,
$$
where $\Vol(h)$ is the volume of the polyhedron $Q(h)$.

\begin{lem}
\label{lem:DVol}
The function $\Vol$ is continuously differentiable on $U$ with
\begin{equation}
\label{eqn:DVol}
\frac{\partial \Vol}{\partial h_i} = A_i,
\end{equation}
where $A_i(h)$ is the area of the $i$-th face of the polyhedron $Q(h)$.
\end{lem}
\begin{proof}
The equation \eqref{eqn:DVol} is geometrically obvious: as we shift the plane of the $i$-th face by
$\epsilon$, we glue to (or cut from) $Q$ a convex slice of thickness $\epsilon$. One side
of the slice has area $A_i$, the other side has area $A_i + O(\epsilon)$. Hense the volume of the
slice is $\epsilon A_i + o(\epsilon)$, and \eqref{eqn:DVol} follows.

Since its partial derivatives $A_i$ are continuous, the function $\Vol$ is continuously
differentiable.
\end{proof}

Equation \eqref{eqn:DVol} implies that the Jacobian of the map $(h_i) \mapsto (A_i)$ equals the
matrix of the second differential of the function $\Vol$. Here by the second differential we mean a
symmetric bilinear form
$$
D^2 \Vol (u,v) := \left. \frac{d}{dt} \right|_{t=0} D_{h+tv} \Vol(u) = \sum_{i,j} \frac{\partial^2
\Vol}{\partial h_i \partial h_j} u_i v_j.
$$
This yields the following reformulation of Theorem \ref{thm:MinkRigDiscr}.

\begin{quote}
\emph{The second differential of the volume at $h=h^0$has corank $3$:
$$
\dim \ker (D_{h^0}^2 \Vol) = 3.
$$
}
\end{quote}

In order to prove Theorem \ref{thm:MinkRigDiscr}, we will compute the second variation of the volume
in two different ways and compare the formulas obtained. Here the \emph{second variation}
$\Vol^{\bigdot\bigdot}$ is the quadratic form associated to $D^2\Vol$:
$$
\Vol^{\bigdot\bigdot} := D^2 \Vol (\dot h, \dot h).
$$
Equivalently,
$$
\Vol^{\bigdot\bigdot} = \left. \frac{d^2}{dt^2} \right|_{t=0} \Vol(h+t\dot h).
$$

\subsection{Orthoscheme decomposition}
\label{subsec:Ortho}
Recall that $F_i$ is the face of $Q(h)$ with outer normal $\nu_i$. If $F_i$ and $F_j$ share an
edge, denote this edge by $F_{ij}$.

Let $q_i$ be the foot of the perpendicular dropped from $0 \in \R^3$ to the plane spanned
by $F_i$. For every pair of adjacent faces, drop perpendiculars from the points $q_i$ and
$q_j$ to the line spanned by $F_{ij}$. Their
common foot will be denoted by $q_{ij}=q_{ji}$. Finally, denote by $q_{ijk}$ the vertex of $Q(h)$
where faces $F_i$, $F_j$, and $F_k$ meet.

\begin{figure}[ht]
\begin{center}
\input{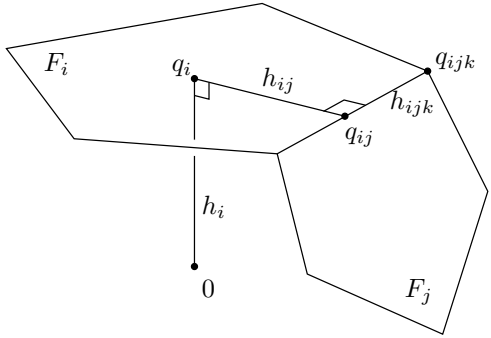}
\end{center}
\caption{Definition of $h_{ij}$ and $h_{ijk}$.}
\label{fig:hijk}
\end{figure}

Denote by $h_{ij}$ the signed length of the segment $q_iq_{ij}$, the sign being positive if $q_i$
lies on the same side from the edge $F_{ij}$ as the polygon $F_i$. Similarly, let $h_{ijk}$ be the
signed length of the segment $q_{ij}q_{ijk}$. See Figure \ref{fig:hijk}.

From now on, we assume that the polyhedron $Q(h^0)$ is \emph{simple}, i.~e. that at each of its
vertices meet exactly three faces. Then, by choosing a neighborhood $U$ of $h^0$ appropriately
small, we can ensure that all polyhedra $Q(h)$ with $h \in U$ are combinatorially isomorphic (under
identification of faces with equal outward normals). The case of a non-simple $Q(h^0)$ is a bit
subtle, and we explain in Subsection \ref{subsec:NonSimple} how it can be treated.

Under this assumption, the functions
\begin{eqnarray*}
h_{ij} \colon U & \to & \R,\\
h_{ijk} \colon U & \to & \R,
\end{eqnarray*}
are linear. If $F_i \cap F_j = \emptyset$ or $F_i \cap F_j \cap F_k = \emptyset$, then the
corresponding functions $h_{ij}$ and $h_{ijk}$ are not defined or can be put identically zero.

\begin{lem}
For all $h \in U$, we have 
\begin{equation}
\label{eqn:VolHHH}
\Vol(h) = \frac16 \sum_{i,j,k} h_i h_{ij} h_{ijk}.
\end{equation}
\end{lem}
\begin{proof}
Denote by $\ell_{ij}$ the length of the edge $F_{ij}$. Then we have
\begin{equation}
\label{eqn:Ell}
\ell_{ij} = h_{ijk} + h_{ijl},
\end{equation}
where $q_{ijk}$ and $q_{ijl}$ are the endpoints of $F_{ij}$. By substituting this into
\begin{equation}
\label{eqn:AHEll}
A_i = \frac12 \sum_j h_{ij} \ell_{ij}
\end{equation}
and substituting the result into
\begin{equation}
\label{eqn:VolHA}
\Vol(h) = \frac13 \sum_i h_i A_i,
\end{equation}
we obtain \eqref{eqn:VolHHH}.

Alternatively, the right hand side of \eqref{eqn:VolHHH} can be seen as the sum of signed volumes
of the orthoschemes $0q_iq_{ij}q_{ijk}$.
\end{proof}

Note that since $h_{ij}$ and $h_{ijk}$ are linear functions of $h$, formula \eqref{eqn:VolHHH}
expresses $\Vol \colon U \to \R$ as a third degree homogeneous polynomial in variables $(h_i)$.

\subsection{First and second variations of the volume and of the face areas}
We start with an auxiliary lemma.
\begin{lem}
\label{lem:Basic}
For an arbitrary variation $\dot h \in \R^n$ we have
\begin{equation}
\label{eqn:hijDot}
\dot h_i h_{ij} + \dot h_j h_{ji} = h_i \dot h_{ij} + h_j \dot h_{ji},
\end{equation}
\begin{equation}
\label{eqn:hijkDot}
\dot h_{ij} h_{ijk} + \dot h_{ik} h_{ikj} = h_{ij} \dot h_{ijk} + h_{ik} \dot h_{ikj}.
\end{equation}
Here $\dot h_{ij}$ and $\dot h_{ijk}$ denote derivatives in the direction of $\dot h$.
\end{lem}
\begin{proof}
Consider the quadrilateral $0q_iq_{ij}q_j$, see Figure \ref{fig:OOij}. The angle $\phi_{ij}$ at $0$
equals the angle between the normals $\nu_i$ and $\nu_j$, and is therefore constant. On the other
hand, we have
\begin{equation}
\label{eqn:Arctan}
\phi_{ij} = \arctan \frac{h_{ij}}{h_i} + \arctan \frac{h_{ji}}{h_j}.
\end{equation}
By differentiating in the direction of $\dot h$, we obtain
\begin{eqnarray*}
0 & = & \frac{1}{1+\frac{h_{ij}^2}{h_i^2}} \frac{\dot h_{ij} h_i - h_{ij} \dot h_i}{h_i^2} +
\frac{1}{1+\frac{h_{ji}^2}{h_j^2}} \frac{\dot h_{ji} h_j - h_{ji} \dot h_j}{h_j^2}\\
& = & \frac{\dot h_{ij} h_i - h_{ij} \dot h_i}{h_i^2 + h_{ij}^2} + \frac{\dot h_{ji} h_j - h_{ji}
\dot h_j}{h_j^2 + h_{ji}^2}.
\end{eqnarray*}
Observing that $h_i^2 + h_{ij}^2 = h_j^2 + h_{ji}^2$, we arrive at \eqref{eqn:hijDot}.

\begin{figure}[ht]
\begin{center}
\input{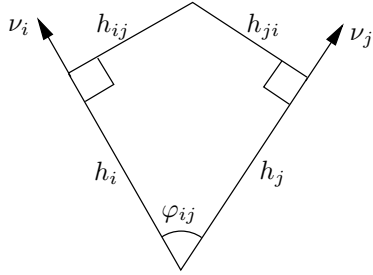}
\end{center}
\caption{The quadrilateral $0q_iq_{ij}q_j$.}
\label{fig:OOij}
\end{figure}

If some of the lengths $h_i$, $h_j$, $h_{ij}$, $h_{ji}$ are negative, the formula \eqref{eqn:Arctan}
still remains valid. Care has to be taken when $h_i = 0$ or $h_j = 0$. One can avoid these
difficulties by noticing that both sides of \eqref{eqn:hijDot} are linear functions of $h_i$ and
$h_j$; so in order to prove \eqref{eqn:hijDot} for all $h$ it suffices to check it for all $h$ in
some open set. By letting $h$ vary over a small neighborhood of the point $(1,1,\ldots,1)$, we
ensure that all segments on Figure \ref{fig:OOij} have positive length. (It does not matter whether
$(1,1,\ldots,1) \in U$, as the functions $(h_{ij})$ can be extended to $\R^n$ by linearity.)

Equation \eqref{eqn:hijkDot} is proved in the same way. For the last part of the argument, one
should note that $(h_{ijk})$ are linear functions of $(h_{ij})$.

Alternatively, the lemma can be proved by using explicit formulas expressing $h_{ij}$ in terms of
$h_i$ and $h_j$, see Figure \ref{fig:Quad}.

\begin{figure}[ht]
\begin{center}
\input{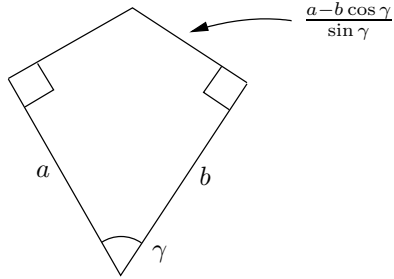}
\end{center}
\caption{Consider orthogonal projection on the side $a$.}
\label{fig:Quad}
\end{figure}
\end{proof}

\begin{lem}
\label{lem:AreaFirst}
The first variations of the face areas are given by
\begin{equation}
\label{eqn:AreaFirst}
\dot A_i = \sum_{j,k} \dot h_{ij} h_{ijk} = \sum_{j,k} h_{ij} \dot h_{ijk}.
\end{equation}
\end{lem}
\begin{proof}
Equations \eqref{eqn:AHEll} and \eqref{eqn:Ell} imply
\begin{eqnarray*}
\dot A_i & = & \frac12 \sum_{j,k} (\dot h_{ij} h_{ijk} + h_{ij} \dot h_{ijk})\\
& = & \frac12 \sum_{\{j,k\}} (\dot h_{ij} h_{ijk} + \dot h_{ik} h_{ikj}) + \frac12 \sum_{\{j,k\}}
(h_{ij} \dot h_{ijk} + h_{ik} \dot h_{ikj}),
\end{eqnarray*}
where $\sum_{\{j,k\}}$ denotes the sum over unordered pairs of $j$ and $k$. By using
\eqref{eqn:hijkDot} and converting back to the sum over ordered pairs, we obtain~
\eqref{eqn:AreaFirst}.
\end{proof}

\begin{lem}
\label{lem:VolFirst}
The first variation of the volume is given by any of the following formulas.
\begin{equation}
\label{eqn:VolFirstA}
\Vol^\bigdot = \sum_i \dot h_i A_i
\end{equation}
\begin{equation}
\label{eqn:VolFirstB}
\Vol^\bigdot = \frac12 \sum_i h_i \dot A_i.
\end{equation}

\end{lem}
\begin{proof}
From \eqref{eqn:VolHA} we have
\begin{equation}
\label{eqn:VolDot}
\Vol^\bigdot = \frac13 \sum_i \dot h_i A_i + \frac13 \sum_i h_i \dot A_i.
\end{equation}
By using \eqref{eqn:AHEll}, regrouping, then using \eqref{eqn:hijDot}, regrouping back, and finally
applying Lemma \ref{lem:AreaFirst}, we obtain
\begin{eqnarray*}
2\sum_i \dot h_i A_i & = & \sum_i \dot h_i \sum_j h_{ij} \ell_{ij} = \sum_{i,j,k} \dot h_i h_{ij}
h_{ijk} = \sum_{\{i,j\},k} (\dot h_i h_{ij} + \dot h_j h_{ji}) h_{ijk}\\
& = & \sum_{\{i,j\},k} (h_i \dot h_{ij} + h_j \dot h_{ji}) h_{ijk} = \sum_{i,j,k} h_i \dot h_{ij}
h_{ijk} = \sum_i h_i \sum_{j,k} \dot h_{ij} h_{ijk}\\
& = & \sum_i h_i \dot A_i.
\end{eqnarray*}
Substituting this in \eqref{eqn:VolDot} yields \eqref{eqn:VolFirstA} and \eqref{eqn:VolFirstB}.
\end{proof}

\begin{rem}
\label{rem:DArea}
Note that \eqref{eqn:VolFirstA} is equivalent to \eqref{eqn:DVol}. Similarly, the first equation in
Lemma \ref{lem:AreaFirst} says that $\frac{\partial A_i}{\partial h_{ij}} = \ell_{ij}$, which is geometrically obvious. This provides an alternative, more geometric approach to Lemmas
\ref{lem:VolFirst} and \ref{lem:AreaFirst}.
\end{rem}

\begin{lem}
The second variations of the face areas are given by
$$
\ddot A_i = \sum_{j,k} \dot h_{ij} \dot h_{ijk}.
$$
\end{lem}
\begin{proof}
This follows from Lemma \ref{lem:AreaFirst} since $\ddot h_{ij} = \ddot h_{ijk} = 0$ because of
their linearity in $h$.
\end{proof}

\begin{lem}
\label{lem:VolSecond}
The second variation of the volume is given by any of the following two formulas.
\begin{equation}
\label{eqn:VolSecondA}
\Vol^{\bigdot\bigdot} = \sum_i \dot h_i \dot A_i.
\end{equation}
\begin{equation}
\label{eqn:VolSecondB}
\Vol^{\bigdot\bigdot} = \sum_i h_i \ddot A_i.
\end{equation}
\end{lem}
\begin{proof}
By differentiating \eqref{eqn:VolFirstA} and taking into account $\ddot h_i = 0$, we
obtain~\eqref{eqn:VolSecondA}. Differentiating \eqref{eqn:VolFirstB} yields
$$
\Vol^{\bigdot\bigdot} = \frac12 \sum_i \dot h_i \dot A_i + \frac12 \sum_i h_i \ddot A_i
$$
Combining this with \eqref{eqn:VolSecondA} yields \eqref{eqn:VolSecondB}.
\end{proof}

\begin{rem}
\label{rem:SymmProp}
Another way to write the first and the second variations of $\Vol$ is
$$
\Vol^\bigdot = \frac12 \sum_{i,j,k} \dot h_i h_{ij} h_{ijk} = \frac12 \sum_{i,j,k} h_i \dot h_{ij}
h_{ijk} = \frac12 \sum_{i,j,k} h_i h_{ij} \dot h_{ijk},
$$
$$
\Vol^{\bigdot\bigdot} = \sum_{i,j,k} \dot h_i \dot h_{ij} h_{ijk} = \sum_{i,j,k} h_i \dot h_{ij}
\dot h_{ijk} = \sum_{i,j,k} \dot h_i h_{ij} \dot h_{ijk}.
$$
\end{rem}

\subsection{Proof of Theorem \ref{thm:MinkRigDiscr}}
\label{subsec:ProofMink}
The proof of Theorem \ref{thm:MinkRigDiscr} is based on the following key lemma.
\begin{lem}
\label{lem:AreaVariationNegative}
Let $\dot h$ be such that $\dot A_i = 0$ for some $i$. Then we have
\begin{equation}
\label{eqn:DDotANegative}
\ddot A_i \le 0.
\end{equation}
Besides, equality in \eqref{eqn:DDotANegative} holds only if $\dot \ell_{ij} = 0$ for all edges
bounding the face~$F_i$.
\end{lem}
This lemma will be proved later in this subsection.

\begin{proof}[Proof of Theorem \ref{thm:MinkRigDiscr}]
We prove Theorem \ref{thm:MinkRigDiscr} in a reformulation given in Subsection \ref{subsec:Theorem}:
if $\dot h$ is such that $\dot A_i = 0$ for all $i$, then $\dot h_i = \langle a, \nu_i \rangle$ for some $a \in \R^3$.

If $\dot A_i = 0$ for all $i$, then \eqref{eqn:VolSecondA} implies
\begin{equation}
\label{eqn:DDotVol0}
\Vol^{\bigdot\bigdot} = \sum_i \dot h_i \dot A_i = 0.
\end{equation}
On the other hand, by Lemma \ref{lem:AreaVariationNegative} we have $\ddot A_i \le 0$. Besides,
without loss of generality we can assume that $h_i > 0$ for all $i$ (just choose the origin $0 \in
\R^3$ inside $Q$). Hence \eqref{eqn:VolSecondB} implies
$$
\Vol^{\bigdot\bigdot} = \sum_i h_i \ddot A_i \le 0.
$$
By comparing this with \eqref{eqn:DDotVol0}, we deduce that $\ddot A_i = 0$ for all $i$. Then, by
the second part of Lemma \ref{lem:AreaVariationNegative}, we have $\dot \ell_{ij} = 0$ for all
edges of $Q$. It follows easily that $\dot h$ translates the polyhedron as a rigid body, thus $\dot h_i = \langle a, \nu_i \rangle$ for some $a \in \R^n$. Theorem \ref{thm:MinkRigDiscr} is proved.
\end{proof}

Fix an index $i$ and consider the function $A_i \colon U \to \R$. It can be written as a function of $(h_{ij})$. Here $j$ varies over all faces of $Q$ adjacent to the $i$-th face. Without loss of
generality assume $j = 1, \ldots, m$ and introduce new variables
$$
g_j = h_{ij},\quad j = 1, \ldots, m.
$$
Formally speaking, we consider a linear map
\begin{eqnarray}
\Phi \colon \R^n & \to & \R^m, \label{eqn:Phi}\\
(h_i) & \mapsto & (g_j) = (h_{ij}). \nonumber
\end{eqnarray}
The map $A_i$ factors through $\Phi$ and induces a map
$$
A \colon V \to \R,
$$
where $V = \Phi(U) \subset \R^m$ is a neighborhood of $g^0 = \Phi(h^0)$. We have
$$
A = \frac12 \sum_{j=1}^m g_j \ell_j,
$$
where $\ell_j = \ell_{ij}$ are linear functions of $g$.

By our construction, $g_j = h_{ij}$ are the support parameters of the face $F_i$ with respect to the projection of $0$ to $\span(F_i)$. We have
\begin{equation}
\label{eqn:DADA}
\frac{\partial A}{\partial g_j} = \ell_j, \mbox{ hence } \frac{\partial^2 A}{\partial g_j \partial g_k} =
\frac{\partial \ell_j}{\partial g_k}.
\end{equation}

\begin{lem}
\label{lem:SignArea}
The second differential $D^2 A$ of the area function on the space of convex $m$--gons with fixed
edge directions has signature $(+_1,0_2,-_{m-3})$. Besides, $D^2 A$ at a point $g$ takes a positive
value on the vector $g$.
\end{lem}
\begin{proof}
First, let us show that
\begin{equation}
\label{eqn:DimKerA}
\dim \ker D^2 A = 2.
\end{equation}
By \eqref{eqn:DADA}, $D^2 A$ coincides with the Jacobian of the map $g \mapsto \ell$. Thus
\begin{equation}
\label{eqn:KerD2A}
\dot g \in \ker D^2 A \ \Longleftrightarrow\ \dot\ell = 0.
\end{equation}
But $\dot\ell = 0$ implies that $\dot g$ is induced by a parallel translation of the polygon. Such
deformations form a 2--dimensional space, and \eqref{eqn:DimKerA} follows.

Consider the space $\Q_m$ of all convex $m$-gons. A convex polygon is determined by its edge normals $(\mu_j)$ and support parameters $(g_j)$. The numbers $(g_j)$ must satisfy a system of linear inequalities (with coefficients depending on $(\mu_j)$) expressing the fact that all edge lengths are non-negative. For every collection $(\mu_j)$ of $m$ different unit vectors positively spanning $\R^2$, the collection $(g_j = 1)$ satisfies this system, as $(\mu, g)$ corresponds to a
circumscribed polygon. Thus $\Q_m$ retracts to the space of configurations $(\mu_j)$ of $m$ unit vectors positively spanning $\R^2$. It follows that $\Q_m$ is connected.

Due to \eqref{eqn:DimKerA}, the rank of $D^2 A$ is constant over the space $\Q_m$. As $\Q_m$ is
connected, it follows that the signature of $D^2 A$ is also constant. Thus in order to determine the signature at the point we need, it suffices to compute it at a point we like.

Let us compute the matrix of $D^2 A$. Let $\alpha_{j,j+1}, j=1, \ldots, m,$ be the exterior angle
between $j$-th and $(j+1)$-st side ($j+1$ taken modulo $m$). Then the formula on Figure \ref{fig:Quad} implies
\begin{equation}
\label{eqn:DEllDG1}
\frac{\partial \ell_j}{\partial g_j} = - (\cot \alpha_{j-1,j} + \cot \alpha_{j,j+1}),
\end{equation}
\begin{equation}
\label{eqn:DEllDG2}
\frac{\partial \ell_j}{\partial g_{j+1}} = \frac{\partial \ell_{j+1}}{\partial g_j} = \frac{1}{\sin
\alpha_{j,j+1}}.
\end{equation}
Note that the matrix of $D^2 A$ does not depend on $g$. This is because $A$ is a homogeneous
polynomial of second degree in $(g_j)$. By the same reason we have
$$
D^2 A(g, g) = 2A > 0,
$$
which proves the second statement in the lemma.

To determine the signature of $D^2A$, put $\alpha_{j,j+1} = \frac{2\pi}{m}$. Then we have
$$
D^2 A = \frac{1}{\sin\frac{2\pi}{m}}
\begin{pmatrix}
-2\cos\frac{2\pi}{m} & 1 & 0 & \ldots & 1\\
1 & -2\cos\frac{2\pi}{m} & 1 & \ldots & 0\\
0 & 1 & -2\cos\frac{2\pi}{m} & \ldots & 0\\
\vdots & \vdots & \vdots & \ddots & \vdots\\
1 & 0 & 0 & \ldots & -2\cos\frac{2\pi}{m}
\end{pmatrix}.
$$ 
The spectrum of this matrix is $\left\{\frac{2(\cos\frac{2\pi k}{m} -
\cos\frac{2\pi}{m})}{\sin\frac{2\pi}{m}}, k = 1, \ldots, m\right\}$, which contains exactly one
positive (for $k = 
m$) and two zero ($k=1$ and $k=m-1$) eigenvalues. The lemma is proved.
\end{proof}

\begin{proof}[Proof of Lemma \ref{lem:AreaVariationNegative}]
Change the variables for the function $A_i$ with the help of the map \eqref{eqn:Phi}. Let $\dot g =
{\mathrm d}\Phi(\dot h)$ be the variation of $g$ induced by $\dot h$. By assumption, we have
$$
\sum_j \dot g_j \ell_j = \dot A = 0.
$$
On the other hand, since $(\ell_j)$ are linear functions of $g$, we have
$$
\ell_j = \sum_k \frac{\partial \ell_j}{\partial g_k} g_k,
$$
which results in
$$
D^2A(\dot g, g) = \sum_{j,k} \frac{\partial \ell_j}{\partial g_k} \dot g_j g_k = \sum_j \dot g_j
\ell_j = 0.
$$
That is, the vectors $\dot g, g \in \R^m$ are mutually orthogonal with respect to the symmetric
bilinear form $D^2A$. By Lemma \ref{lem:SignArea}, $g$ is a positive vector for $D^2A$, and $D^2A$
is negative semidefinite on the orthogonal complement to $g$. Thus we have
$$
\ddot A = D^2A(\dot g, \dot g) \le 0.
$$
If $\ddot A = 0$, then $\dot g \in \ker D^2A$, which by \eqref{eqn:KerD2A} implies $\dot \ell_j = 0$
for all $j$. The lemma is proved.
\end{proof}

\begin{rem}
In general, the second differential does not behave well under change of variables. However, it does in our case, because the map \eqref{eqn:Phi} is linear. Therefore the second variation $\ddot A_i$ in Lemma \ref{lem:AreaVariationNegative} is associated with the second differential $D^2 A$ from Lemma \ref{lem:SignArea}:
$$
\ddot A_i = D^2A({\mathrm d}\Phi(\dot h), {\mathrm d}\Phi(\dot h)) = D^2A(\dot g, \dot g),
$$
which was implicitly used in the above proof.
\end{rem}

\subsection{The case of a non-simple polyhedron $Q$}
\label{subsec:NonSimple}
In the course of our proof we assumed that the polyhedron $Q$ is simple, that is each of its
vertices belongs to exactly three faces. This assumption implies that the combinatorics of $Q$ is
preserved when its support parameters vary slightly.

If the polyhedron $Q = Q(h^0)$ is non-simple, then $Q(h)$ may be combinatorially different from $Q$
even for $h$ close to $h^0$. A neighborhood $U$ of $h^0 \in \R^n$ is subdivided into cells
$(U^\Delta)$ with $\Delta$ indexing simple combinatorial types of perturbed polyhedra. On each cell,
the function $\Vol$ is a third degree homogeneous polynomial in $(h_i)$:
$$
\Vol(h) = V^\Delta(h) \mbox{ for } h \in U^\Delta,
$$
so that function $\Vol$ is piecewise polynomial in a neighborhood of $h^0$. It can be shown that
$$
\Vol \in C^2(U), \quad A_i \in C^1(U).
$$
But $A_i$ may fail to be $C^2$, which is bad because second differentials of face areas play a key
role in Subsection \ref{subsec:ProofMink}.

This problem can be resolved as follows. Choose any simple combinatorics $\Delta$ from a neigborhood of $h^0$, and work with $V^\Delta$ and $A_i^\Delta$ instead of $\Vol$ and $A_i$. Geometrically this means that we view our non-simple polyhedron $Q(h^0)$ as a member of a family of simple polyhedra, with some edge lengths vanishing; by varying $h^0$ in $U$, we allow negative edge
lengths to appear. One can check that all arguments in Subsection \ref{subsec:ProofMink} go through.

This modification proves Theorem \ref{thm:MinkRigDiscr} for non-simple polyhedra.

\section{Metric rigidity of convex polyhedra}
\label{sec:MetrRig}
\subsection{The theorem}
\label{subsec:ThmMetric}
Let $P \subset \R^3$ be a convex polyhedron with vertices $p_1, \ldots, p_n$. Assume that $P$ is
simplicial, that is all of its faces are triangles. Subsection \ref{subsec:NonSimplicial} explains
how our arguments change in the non-simplicial case.

\begin{dfn}
\label{dfn:InfDeform}
An \emph{infinitesimal deformation} of $P$ is a collection of vectors $q_i \in \R^3, i=1, \ldots,
n$. Each $q_i$ is thought of as a vector applied at the point~$p_i$.

An infinitesimal deformation $(q_i)$ of a polyhedron $P$ is called \emph{isometric}, if
\begin{equation}
\label{eqn:LDT}
\left. \frac{d}{dt} \right|_{t=0} \| (p_i+tq_i) - (p_j+tq_j) \| = 0,
\end{equation}
for all edges $p_ip_j$ of $P$. In other words, if lengths of all edges don't change in the first
order as vertices $(p_i)$ move linearly with velocities $(q_i)$.
\end{dfn}

A simple computation shows that condition \eqref{eqn:LDT} is equivalent to
$$
\langle p_i - p_j, q_i - q_j \rangle = 0.
$$

\begin{dfn}
An infinitesimal deformation $(q_i)$ of $P$ is called \emph{trivial}, if the map $p_i \mapsto q_i$
is the restriction of an infinitesimal isometry of $\R^3$:
$$
q_i = A p_i + b, \quad A \in \mathfrak{so}(3), b \in \R^3.
$$
\end{dfn}

An infinitesimal isometry of $\R^3$ preserves in the first order distances between \emph{all} pairs
of points; therefore every trivial infinitesimal deformation of $P$ is isometric.

\begin{dfn}
A polyhedron $P$ is called \emph{infinitesimally rigid}, if every isometric infinitesimal
deformation of $P$ is trivial.
\end{dfn}

\begin{thm}[Legendre-Cauchy-Dehn]
\label{thm:PolyhRig}
Every convex polyhedron in $\R^3$ is infinitesimally rigid.
\end{thm}

This is usually referred to as Dehn's theorem. Let us explain why we prefer a different attribution.
Cauchy \cite{Cau13} proved a global rigidity statement: two convex polyhedra with the same
combinatorics and pairwise isometric faces are congruent. His proof was based on ideas presented by
Legendre in \cite[note XII, pp. 321--334]{Leg94}. Dehn was the first to state and prove the
infinitesimal rigidity theorem in \cite{Dehn16}. However, in the footnote on the first page of
\cite{Weyl17}, Weyl remarks that the argument in \cite{Cau13} carries over word by word to yield a
proof of Theorem \ref{thm:PolyhRig} and reproaches Dehn for not citing Cauchy. Note also that the
Cauchy's ``arm lemma'' is more immediate in the infinitesimal context.

\subsection{A reformulation}
\label{subsec:WarpPolyh}
Without loss of generality we may assume that the coordinate origin $0\in\R^3$ lies in the interior
of $P$. For every vertex $p_i$ of $P$, put
$$
r_i^0 = \|p_i\|.
$$
The length of an edge $p_ip_j$ will be denoted by
$$
\ell_{ij} = \|p_i - p_j\|.
$$
Triangles $0p_ip_j$, where $p_ip_j$ ranges over all edges of $P$, cut the polyhedron $P$ into
triangular pyramids. These pyramids have a common apex at $0$, their bases are faces of $P$. Let $r
\in \R^n$ be a point close to $r^0 = (r_i^0)_{i=1}^n$. Change lateral edge lengths of the pyramids
from $r_i^0$ to $r_i$ while keeping base edge lengths $\ell_{ij}$ fixed. A metric space glued from
the new collection of pyramids (by the old gluing rules) is called a \emph{warped polyhedron}. A
warped polyhedron is not embeddable in $\R^3$ in general, because the sum $\omega_i$ of all dihedral
angles at an edge $0p_i$ may be different from $2\pi$, for some $i$. Denote by
$$
\kappa_i = 2\pi - \omega_i
$$
the \emph{curvature} of the warped polyhedron at the edge $0p_i$. This yields a $C^\infty$-map
\begin{eqnarray*}
U & \to & \R^n,\\
r & \mapsto & \kappa,
\end{eqnarray*}
where $U \subset \R^n$ is a sufficiently small neighborhood of $r^0$.

\begin{lem}
\label{lem:Reformulation}
A convex polyhedron $P$ is infinitesimally rigid if and only if
$$
\dim\ker \left.\left(\frac{\partial \kappa_i}{\partial r_j}\right)\right|_{r=r^0} = 3.
$$
\end{lem}
\begin{proof}
We will establish a correspondence between isometric infinitesimal deformations of $P$ and elements
of the kernel of $J_r^\kappa = \left(\frac{\partial \kappa_i}{\partial r_j}\right)$.

Let $q$ be an infinitesimal isometric deformation of $P$ and let $p_i^t = p_i + tq_i$ be the
corresponding linear motions of the vertices. Put
\begin{equation}
\label{eqn:Deform}
r_i^t = \|p_i^t\|,\quad \ell_{ij}^t = \|p_i^t - p_j^t\|.
\end{equation}
This defines a family of warped polyhedra, this time with non-constant metric on the boundary. By
Definition \ref{dfn:InfDeform} and since the polyhedron remains embedded in $\R^3$, we have
$$
\dot \ell_{ij} := \left. \frac{d}{dt} \right|_{t=0} \ell_{ij}^t = 0, \quad \kappa_i^t = 0
$$
for all $i,j$. Modify the deformation \eqref{eqn:Deform} by putting $\ell_{ij}^t = \ell_{ij}$.
Although the functions $\kappa_i^t$ might be not identically zero anymore, their time derivatives at
$t=0$ will not change:
$$
\ell_{ij}^t = \ell_{ij} \Rightarrow \dot \kappa_i = 0.
$$
It follows that $\dot r \in \ker J_r^\kappa$.

In particular, if $q$ is an infinitesimal rotation around an axis through the origin, then we have
$r_i^t = r_i$, and thus $\dot r = 0$. If $q$ is a parallel translation, then $\dot r \ne 0$ is a
non-trivial element in $\ker J_r^\kappa$. One can show that translations in linearly independent
directions in $\R^3$ produce linearly independent vectors $\dot r \in \R^n$. Hence trivial
infinitesimal deformations of $P$ give rise to a 3--dimensional subspace of $\ker J_r^\kappa$.

In the opposite direction, let us associate to $\dot r \in \ker J_r^\kappa$ an isometric
infinitesimal deformation of $P$. Let $p_1p_2p_3$ be a face of $P$. Choose a vector $q_1$ collinear
with $0p_1$ so that when $p_1$ moves with the velocity $q_1$, its distance from $0$ changes with the
speed $\dot r_1$:
$$
q_1 = \dot r_1 \frac{p_1}{\|p_1\|}.
$$
Choose $q_2$ collinear with the plane $0p_1p_2$ so that to satisfy the conditions on variations of
$r_2$ and $\ell_{12}$:
$$
\left\langle q_2, \frac{p_2}{\|p_2\|} \right\rangle = \dot r_2, \quad \langle p_1-p_2, q_1-q_2
\rangle = 0.
$$
And choose $q_3$ so that to satisfy the conditions on variations of $r_3$, $\ell_{13}$,
and~$\ell_{23}$:
\begin{equation}
\label{eqn:DetermineQ}
\left\langle q_3, \frac{p_3}{\|p_3\|} \right\rangle = \dot r_3, \quad \langle p_1-p_3, q_1-q_3
\rangle = 0, \quad \langle p_2-p_3, q_2-q_3 \rangle = 0.
\end{equation}
Now proceed to an adjacent face $p_1p_3p_4$ and determine $q_4$ by conditions similar to
\eqref{eqn:DetermineQ}. Continue the face path $(p_1p_ip_{i+1})$ around the vertex $p_1$ until it
closes at $p_1p_kp_2$. On the vector $q_k$ there will be four conditions, the fourth being the
stability of $\ell_{k2}$. It is easy to see that this condition follows from the first three due to
$$
\dot\kappa_1 := \sum_i \frac{\partial \kappa_1}{\partial r_i} \dot r_i = 0.
$$
so that $q_k$ is well-defined. In a similar way we can find $q_i$ for all $i$, all closing
conditions being satisfied due to $\dot \kappa = 0$.

Note that we made some voluntary choices at the beginning, with $q_1$ and $q_2$, but starting from
$q_3$ everything was forced. It follows that the infinitesimal deformation $q$ is determined by
$\dot r$ uniquely up to an infinitesimal rotation.

As we observed, trivial infinitesimal deformations of $P$ generate a 3--dimensional subspace of
$\ker J_r^\kappa$. If there is a non-trivial isometric infinitesimal deformation, then the
corresponding $\dot r \in \ker J_r^\kappa$ lies outside this subspace, as the map $q \mapsto \dot r$
identifies only deformations that differ by an infinitesimal rotation. Thus $\dim \ker J_r^\kappa =
3$ implies infinitesimal rigidity. Vice versa, if $\dim \ker J_r^\kappa > 3$, then there is an $\dot
r \in \ker J_r^\kappa$ not coming from a trivial infinitesimal deformation. Thus it determines a
non-trivial isometric infinitesimal deformation $q$. That is, infinitesimal rigidity implies $\dim
\ker J_r^\kappa = 3$.
\end{proof}

\subsection{The discrete Hilbert-Einstein functional}
\label{subsec:DHE}
In Section \ref{subsec:WarpPolyh} we introduced warped polyhedra that have curvatures $(\kappa_i)$
along their radial edges $0p_i$. Recall that the boundary edge lengths $(\ell_{ij})$ are kept
constant, while the radial edge lengths $r=(r_i)_{i=1}^n$ are allowed to vary in a neighborhood $U$
of a point $r^0 \in \R^n$. Denote by $P(r)$ the warped polyhedron with radial edge lengths $r$. Let
$\lambda_{ij}$ be the exterior dihedral angle at the boundary edge $p_ip_j$. Note that
$(\lambda_{ij})$, as well as $(\kappa_i)$ are functions of $r$.

\begin{dfn}
The \emph{discrete Hilbert-Einstein functional} of a warped polyhedron $P(r)$ is
\begin{equation}
\label{eqn:DHE}
\HE(r) = \sum_i r_i \kappa_i + \sum_{\{i,j\}} \ell_{ij} \lambda_{ij}.
\end{equation}
\end{dfn}
By $\sum_{\{i,j\}}$ we denote the sum over all unordered pairs of $i$ and $j$, so that every
boundary edge is taken once.

\begin{lem}
\label{lem:SFirstVar}
We have
\begin{equation}
\label{eqn:dSr}
\frac{\partial\HE}{\partial r_i} = \kappa_i.
\end{equation}
\end{lem}

Lemma \ref{lem:SFirstVar} can be derived from (and is equivalent to) the Schl\"afli formula for
Euclidean polyhedra. We will prove it by a different method in Subsection \ref{subsec:HEVariations}.

\begin{cor}
\label{cor:CorankDHE}
A convex polyhedron $P$ is infinitesimally rigid if and only if the second differential of the
Hilbert-Einstein functional has corank $3$ at the point $r^0$:
$$
\dim \ker (D^2\HE)|_{r=r^0} = 3.
$$
\end{cor}
\begin{proof}
In view of Lemma \ref{lem:SFirstVar}, the matrix of the second differential of $\HE$ is the Jacobi
matrix of the map $r \mapsto \kappa$. Thus Corollary \ref{cor:CorankDHE} is simply a reformulation
of Lemma \ref{lem:Reformulation}.
\end{proof}

\subsection{The total curvature of a warped spherical polygon}
\label{subsec:TotCurv}
\begin{figure}[ht]
\begin{center}
\input{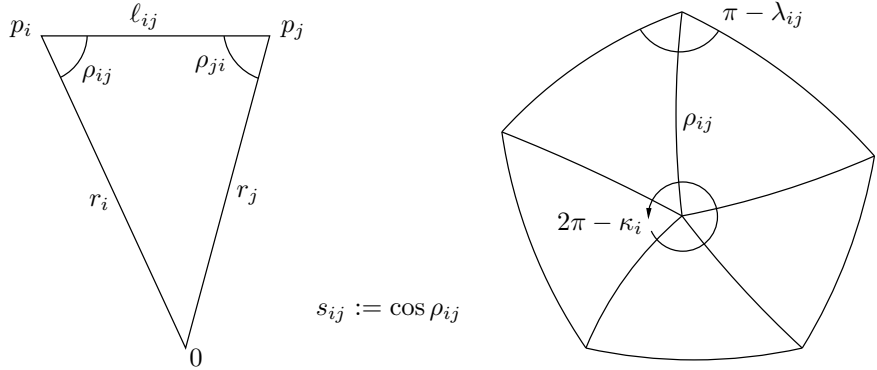}
\end{center}
\caption{The triangle $0p_ip_j$ and the spherical link $S_i$ of the vertex $p_i$.}
\label{fig:SpherLink}
\end{figure}

In a warped polyhedron $P(r)$, consider the spherical link $S_i$ of the vertex $p_i$. By definition,
$S_i$ is a complex of spherical triangles obtained from trihedral angles at $p_i$ by intersecting
them with a unit sphere centered at $p_i$. In our case, these triangles are glued cyclically around
a common vertex, forming a total angle of $2\pi - \kappa_i$. The angles at the boundary vertices of
$S_i$ are equal to the dihedral angles of $P(r)$, that is to $\pi-\lambda_{ij}$. The lengths of
radial edges in $S_i$ are equal to $\rho_{ij}$, where $\rho_{ij}$ is the angle at the vertex $p_i$
in the triangle $0p_ip_j$. See Figure \ref{fig:SpherLink}. Denote
$$
s_{ij} = \cos \rho_{ij}.
$$

A complex of spherical triangles similar to $S_i$ (a set of triangles glued cyclically around a
common vertex) will be called a \emph{warped spherical polygon}. Note that $S_i$ is in addition
\emph{convex}, that is the angles at its boundary vertices are less or equal $\pi$.

\begin{dfn}
The \emph{total curvature} $K_i$ of a warped spherical polygon $S_i$ on Figure \ref{fig:SpherLink},
right, is defined as
$$
K_i = \kappa_i + \sum_j s_{ij} \lambda_{ij}.
$$
\end{dfn}

\begin{lem}
We have
\begin{equation}
\label{eqn:HEK}
\HE = \sum_i r_i K_i.
\end{equation}
\end{lem}
\begin{proof}
The orthogonal projection of $0$ to the line $p_ip_j$ splits the edge $p_ip_j$ in two segments of
lengths $r_i\cos\rho_{ij}$ and $r_j\cos\rho_{ji}$. Thus we have
\begin{equation}
\label{eqn:LRQ}
\ell_{ij} = r_i s_{ij} + r_j s_{ji}.
\end{equation}
(If the projection of $0$ lies outside the edge, the formula remains valid.) By substituting
\eqref{eqn:LRQ} in the definition of $\HE$, we obtain
\begin{eqnarray*}
\HE & = & \sum_i r_i \kappa_i + \sum_{\{i,j\}}(r_is_{ij} + r_js_{ji})\lambda_{ij} = \sum_i r_i
\kappa_i + \sum_{i,j}r_i s_{ij} \lambda_{ij}\\
& = & \sum_i r_i \Big(\kappa_i + \sum_j s_{ij} \lambda_{ij}\Big) = \sum_i r_i K_i,
\end{eqnarray*}
and the lemma is proved.
\end{proof}

\subsection{Variations of $\HE$ and $K_i$}
\label{subsec:HEVariations}
Let $r \in U$, and let $\dot r \in \R^n$ be an arbitrary variation of $r$. Denote by
$$
\HE^\bigdot = \left. \frac{d}{dt} \right|_{t=0} \HE(r+t\dot r)
$$
the derivative of the function $\HE$ in the direction $\dot r$. Similar notations $\dot K_i$, $\dot
s_{ij}$, $\dot \lambda_{ij}$ are used for the directional derivatives of other functions of $r$.

\begin{lem}
\label{lem:KFirstVar}
The first variation of the total curvature of the spherical link $S_i$ of the vertex $p_i$ in $P(r)$
is given by
$$
\dot K_i = \sum_j \dot s_{ij} \lambda_{ij}.
$$
\end{lem}
\begin{proof}
The warped polygon $S_i$ is made of spherical triangles, see Figure \ref{fig:SpherLink}. Consider a
triangle with edge lengths $\rho_{ij}$, $\rho_{ik}$. Its third side has a constant length, as it
equals to the angle $p_jp_ip_k$ on the boundary of the polyhedron $P(r)$. By applying Lemma
\ref{lem:DotABC}, we obtain
\begin{equation}
\label{eqn:DotInS}
\dot \alpha_{ijk} + \dot \beta_{ij} s_{ij} + \dot \gamma_{ik} s_{ik} = 0,
\end{equation}
where $\alpha_{ijk}$ is the angle between the $\rho_{ij}$ and $\rho_{ik}$ sides, and $\beta_{ij}$,
$\gamma_{ik}$ are the two other angles. As we have
$$
\sum_{jk} \alpha_{ijk} = 2\pi - \kappa_i, \quad \beta_{ij} + \gamma_{ij} = \pi - \lambda_{ij},
$$
summing \eqref{eqn:DotInS} over all triangles yields
\begin{equation}
\label{eqn:DotKappa}
\dot \kappa_i + \sum_j s_{ij} \dot \lambda_{ij} = 0.
\end{equation}
It follows that
$$
\dot K_i = \dot \kappa_i + \sum_j (\dot s_{ij} \lambda_{ij} + s_{ij} \dot\lambda_{ij}) = \sum_j \dot
s_{ij} \lambda_{ij},
$$
and the lemma is proved.
\end{proof}

\begin{lem}
\label{lem:RQ}
For every edge $p_ip_j$, we have
$$
\dot r_i s_{ij} + \dot r_j s_{ji} + r_i \dot s_{ij} + r_j \dot s_{ji} = 0.
$$
\end{lem}
\begin{proof}
This follows from \eqref{eqn:LRQ} and $\dot\ell_{ij} = 0$.
\end{proof}

\begin{lem}
\label{lem:HEFirstVar}
The first variation of the Hilbert-Einstein functional of a warped polyhedron $P(r)$ is given by
\begin{equation}
\label{eqn:HEDot}
\HE^\bigdot = \sum_i \dot r_i \kappa_i.
\end{equation}
\end{lem}
\begin{proof}
Follows by differentiating \eqref{eqn:HEK} and applying Lemmas \ref{lem:KFirstVar} and \ref{lem:RQ}:
\begin{eqnarray*}
\HE^\bigdot & = & \sum_i \dot r_i K_i + \sum_i r_i \dot K_i = \sum_i \dot r_i (\kappa_i + \sum_j
s_{ij}\lambda_{ij}) + \sum_i r_i \sum_j \dot s_{ij} \lambda_{ij}\\
& = & \sum_i \dot r_i \kappa_i + \sum_{i,j} (\dot r_i s_{ij} + r_i \dot s_{ij}) \lambda_{ij}\\
& = & \sum_i \dot r_i \kappa_i + \sum_{\{i,j\}} (\dot r_i s_{ij} + \dot r_j s_{ji} + r_i \dot s_{ij}
+ r_j \dot s_{ji}) \lambda_{ij} = \sum_i \dot r_i \kappa_i.
\end{eqnarray*}
\end{proof}

Lemma \ref{lem:SFirstVar} follows from Lemma \ref{lem:HEFirstVar}, as \eqref{eqn:dSr} is just a
reformulation of~\eqref{eqn:HEDot}.

Recall that the \emph{second variation} of $\HE$ is the quadratic form in $(\dot r_i)$ associated to
the second differential of $\HE$:
$$
\HE^{\bigdot\bigdot} = D^2 \HE (\dot r, \dot r) = \sum_{i,j} \frac{\partial^2 \HE}{\partial r_i
\partial r_j} \dot r_i \dot r_j = \left. \frac{d^2}{dt^2} \right|_{t=0} \HE(r + t\dot r).
$$

\begin{lem}
\label{lem:SDDot1}
The second variation of $\HE$ is given by any of the following two formulas.
\begin{equation}
\label{eqn:SDDot1}
\HE^{\bigdot\bigdot} = \sum_i \dot r_i \dot\kappa_i.
\end{equation}
\begin{equation}
\label{eqn:SDDot2}
\HE^{\bigdot\bigdot} = \sum_i r_i \sum_j \dot s_{ij} \dot\lambda_{ij}.
\end{equation}
\end{lem}
\begin{proof}
Equation \eqref{eqn:SDDot1} is a direct consequence of \eqref{eqn:HEDot}. In order to prove
\eqref{eqn:SDDot2}, transform \eqref{eqn:SDDot1} by using \eqref{eqn:DotKappa} and Lemma
\ref{lem:RQ}:
\begin{eqnarray*}
\HE^{\bigdot\bigdot} & = & \sum_i \dot r_i \dot\kappa_i = - \sum_i \dot r_i \sum_j s_{ij} \dot\lambda_{ij} = -
\sum_{\{i,j\}} (\dot r_i s_{ij} + \dot r_j s_{ji}) \dot\lambda_{ij}\\
& = & \sum_{\{i,j\}}(r_i \dot s_{ij} + r_j \dot s_{ji}) \dot\lambda_{ij} = \sum_i r_i \sum_j \dot
s_{ij} \dot\lambda_{ij}.
\end{eqnarray*}
\end{proof}

\subsection{Proof of Theorem \ref{thm:PolyhRig}}
The proof is based on the following lemma.
\begin{lem}
\label{lem:KDDotNegative}
Let $\dot r$ be such that $\dot \kappa_i = 0$ for some $i$. Then we have
\begin{equation}
\label{eqn:SDotLDot}
\sum_j \dot s_{ij} \dot \lambda_{ij} \le 0.
\end{equation}
Besides, equality in \eqref{eqn:SDotLDot} holds only if $\dot \lambda_{ij} = 0$ for all $j$.
\end{lem}
The proof of this lemma requires a detailed study of the second variation of $K_i$ with respect to
the variables $(s_{ij})$. This is done in Subsection \ref{subsec:SecVarK} which ends with the proof
of Lemma \ref{lem:KDDotNegative}. Now we prove Theorem \ref{thm:PolyhRig} assuming the validity of
Lemma \ref{lem:KDDotNegative}.

\begin{proof}[Proof of Theorem \ref{thm:PolyhRig}]
We prove theorem \ref{thm:PolyhRig} in the reformulation obtained in Lemma \ref{lem:Reformulation}.
It suffices to show that if $\dot r$ is such that $\dot \kappa_i = 0$ for all $i$, then $\dot r$ is
induced by a trivial isometric deformation of $P$.

If $\dot \kappa_i = 0$ for all $i$, then \eqref{eqn:SDDot1} implies
$$
\HE^{\bigdot\bigdot} = \sum_i \dot r_i \dot\kappa_i = 0.
$$
On the other hand, equation \eqref{eqn:SDDot2} and Lemma \ref{lem:KDDotNegative} imply
\begin{equation}
\label{eqn:4}
\HE^{\bigdot\bigdot} = \sum_i r_i \sum_j \dot s_{ij} \dot \lambda_{ij} \le 0,
\end{equation}
because $r_i > 0$ for all $i$. Thus we must have
$$
\sum_j \dot s_{ij} \dot \lambda_{ij} = 0
$$
for all $i$. By Lemma \ref{lem:KDDotNegative}, this happens only if $\dot \lambda_{ij} = 0$ for all
edges $p_ip_j$. This means that the isometric infinitesimal deformation that corresponds to $\dot r$
does not change the dihedral anges of the polyhedron $P$ in the first order. It follows that this
infinitesimal deformation is trivial.
\end{proof}

\subsection{The second variation of the total curvature of a warped spherical polygon}
\label{subsec:SecVarK}
In Subsection \ref{subsec:TotCurv}, we defined warped spherical polygons and their total curvature.
Warped spherical polygons appeared as spherical links $S_i$ of vertices of a warped polyhedron
$P(r)$, see Figure \ref{fig:SpherLink}. Consequently, the total curvature of $S_i$ was viewed as a
function of $r \in \R^n$.

In this subsection, we study warped spherical polygons on their own. We preserve the notations on
the right of Figure \ref{fig:SpherLink}, but suppress the index $i$. Without loss of generality,
assume that the index $j$ numbering the radial edges varies from $1$ to $m$. As before, lengths of
boundary edges are assumed constant, so that a warped spherical polygon is determined by $m$
parameters
$$
s_j = \cos \rho_j, \quad j = 1, \ldots, m,
$$
with $s$ varying in a neigborhood $V$ of $s^0 \in \R^m$. The point $s^0$ has the property
$$
\kappa(s^0) = 0,
$$
where $\kappa$ (a former $\kappa_i$) is the angular defect of the interior vertex of the warped
polygon.

Consider the total curvature of a warped spherical polygon as a function of $s$:
\begin{eqnarray}
K \colon V & \to & \R, \nonumber \\
K(s) & = & \kappa + \sum_j s_j \lambda_j. \label{eqn:KInS}
\end{eqnarray}

\begin{lem}
\label{lem:KVarInS}
The first and second variations of $K$ as a function of $s$ are given by
\begin{equation}
\label{eqn:DotKInS}
\dot K = \sum_j \dot s_j \lambda_j,
\end{equation}
\begin{equation}
\label{eqn:DDotKInS}
\ddot K = \sum_j \dot s_j \dot \lambda_j.
\end{equation}
\end{lem}
\begin{proof}
Equation \eqref{eqn:DotKInS} follows from Lemma \ref{lem:KFirstVar}, as the form of the first
variation does not depend on the choice of variables. (Actually, in the proof of Lemma
\ref{lem:KFirstVar} we deal with coordinates $(s_{ij})$, so we implicitly used the invariance of the
first variation.)

Equation \eqref{eqn:DDotKInS} is a direct consequence of \eqref{eqn:DotKInS}. Compare this with
\eqref{eqn:SDDot1}.
\end{proof}

We have abused notation: in \eqref{eqn:DDotKInS}, double dots denote the second variation in the
variables $(s_j)$, while before we used them to denote the second variation in variables $(r_i)$.
Don't ever try to substitute \eqref{eqn:DDotKInS} in \eqref{eqn:SDDot2}!

\begin{lem}
\label{lem:SignD2K}
The second differential $D_{s^0}^2K$ of the function \eqref{eqn:KInS} at the point $s^0$ has
signature $(+_1, 0_2, -_{m-3})$. Besides, the associated quadratic form is positive on the vector
$s^0$:
\begin{equation}
D_{s^0}^2K(s^0, s^0) > 0.
\end{equation}
\end{lem}
\begin{proof}
First let us show that
\begin{equation}
\label{eqn:DimKerK}
\dim \ker D_{s(0)}^2K = 2.
\end{equation}
Lemma \ref{lem:KVarInS} implies
$$
\frac{\partial K}{\partial s_j} = \lambda_j, \mbox{ hence } \frac{\partial^2 K}{\partial s_j
\partial s_k} = \frac{\partial \lambda_j}{\partial s_k},
$$
so that the matrix of $D^2K$ coincides with the Jacobi matrix of the map $s \mapsto \lambda$. Thus
$$
\dot s \in \ker D^2K \Longleftrightarrow \dot \lambda = 0,
$$
that is if and only if the variation $\dot s$ of radial edge lengths induces a zero variation of
angles between boundary edges. It is not hard to show that at $s=s^0$ such variations form a
2--dimensional space, namely they come from moving the interior vertex while the boundary of the
polygon remains fixed.

Let $\P_m$ be the space of all convex spherical $m$--gons with a marked interior point (i.~e. convex
warped spherical $m$--gons with $\kappa = 0$). It is easy to see that $\P_m$ is connected. Every
element of $\P_m$ can be viewed as a point $s^0$ in a space of warped polygons with a fixed boundary
metric, thus has a matrix $D_{s^0}^2 K$ associated with it. By the previous paragraph, the rank of
$D_{s^0}^2 K$ is constant over $\P_m$. Since $\P_m$ is connected, and $D_{s^0}^2 K$ depends
continuously on an element of $\P_m$, the signature is constant as well. Thus in order to compute
the signature of $D_{s^0}^2 K$ for \emph{all} warped spherical $m$--gons with $\kappa = 0$, it
suffices to do this for \emph{one} such $m$--gon.

Let us compute the matrix of $D_{s^0}^2 K$. Let $\alpha_{j, j+1}$ be the angle between two
consecutive radial edges. Then we have by Lemma \ref{lem:DADAlpha}
$$
\frac{\partial \lambda_j}{\partial s_j} = - \frac{\cot \alpha_{j-1,j} + \cot \alpha_{j,j+1}}{\sin^2
\rho_j},
$$
$$
\label{eqn:DLambdaJK}
\frac{\partial \lambda_j}{\partial s_{j+1}} = \frac{\partial \lambda_{j+1}}{\partial s_j} =
\frac{1}{\sin\alpha_{j,j+1} \sin\rho_j \sin\rho_{j+1}}.
$$

We compute the signature of $D_{s^0}^2 K$ for $\alpha_{j, j+1} = \frac{2\pi}{m}$ and $\rho_j = \rho
\in (0,\frac{\pi}{2})$. In this case
$$
D_{s^0}^2 K = \frac{1}{\sin\frac{2\pi}{m} \sin^2 \rho}
\begin{pmatrix}
-2\cos\frac{2\pi}{m} & 1 & 0 & \ldots & 1\\
1 & -2\cos\frac{2\pi}{m} & 1 & \ldots & 0\\
0 & 1 & -2\cos\frac{2\pi}{m} & \ldots & 0\\
\vdots & \vdots & \vdots & \ddots & \vdots\\
1 & 0 & 0 & \ldots & -2\cos\frac{2\pi}{m}
\end{pmatrix}.
$$ 
This matrix has spectrum $\left\{\frac{2(\cos\frac{2\pi k}{m} -
\cos\frac{2\pi}{m})}{\sin\frac{2\pi}{m} \sin^2 \rho}, k = 1, \ldots, m\right\}$, which contains
exactly one positive eigenvalue (for $k=m$). The first part of Lemma \ref{lem:SignD2K}, concerning
the signature, is proved.

The positivity of the form $D_{s^0}^2K$ on the vector $s^0$ follows from Lemma \ref{lem:DA=AreaDual}
below.
\end{proof}

Let $S \subset \Sph^2$ be a convex spherical polygon with a distinguished point $p_0$ in its
interior. The polar dual $S^*$ of $S$ is the intersection of all hemispheres centered at vertices of
$S$. We define the \emph{Euclidean polar dual} $S^*_\Euc$ of $S$ as the projection of $S^*$ from the
center of $\Sph^2$ to the plane tangent to $\Sph^2$ at $p_0$. See Figure \ref{fig:EuclPolarDual}.
Note that if some of the distances $\rho_j$ are bigger than $\frac\pi 2$, then the point $p_0$ lies
outside the spherical and Euclidean polar duals.

\begin{figure}[ht]
\begin{center}
\input{EuclPolDual.tex}
\end{center}
\caption{Fragments of $S$, $S^*$, and $S^*_\Euc$.}
\label{fig:EuclPolarDual}
\end{figure}

\begin{lem}
\label{lem:DA=AreaDual}
Let $S$ be a convex spherical polygon with a marked interior point $p_0$. Let $(s_j^0)$ be distances
from $p_0$ to the vertices of $S$. Consider the space of warped spherical polygons obtained from $S$
by varying $(s_j)$ while keeping the boundary edge lengths fixed. Then we have
$$
D_{s^0}^2K(s^0, s^0) = 2 \operatorname{Area}(S^*_\Euc),
$$
where $S^*_\Euc$ is the Euclidean polar dual to $S$, see the definition before the lemma.
\end{lem}
\begin{proof}
Draw the perpendiculars from $p_0$ to the sides of $S^*$. They cut the spherical polygon $S^*$ into
quadrilaterals, each with a pair of opposite right angles (quadrilaterals may be self-intersecting).
The signed lengths of the perpendiculars equal $\frac{\pi}{2} - \rho_j$, with angles $\alpha_{j,
j+1}$ between them. The Euclidean polygon $S^*_\Euc$ has a similar decomposition, with the same
angles $\alpha_{j, j+1}$ but with perpendiculars of lengths $\cot \rho_j$, see Figure
\ref{fig:EuclPolarDual}. Formula on Figure \ref{fig:Quad} gives us the lengths of the other two sides of a $(j,j+1)$-quadrilateral, so that we can compute its area:
$$
\frac12 \left( \cot\rho_j \frac{\cot\rho_{j+1} - \cot\rho_j \cos\alpha_{j,j+1}}{\sin \alpha_{j,j+1}}
+ \cot\rho_{j+1} \frac{\cot\rho_j - \cot\rho_{j+1} \cos\alpha_{j,j+1}}{\sin \alpha_{j,j+1}} \right).
$$
By summing over $j$ and performing simple transformations, we obtain
$$
\operatorname{Area}(S^*_\Euc) = \frac{1}{2}\sum_{j,k} \frac{\partial \lambda_j}{\partial s_k} s_j^0
s_k^0 = \frac12 D_{s^0}^2K(s^0,s^0),
$$
and the lemma is proved.
\end{proof}

\begin{proof}[Proof of Lemma \ref{lem:KDDotNegative}]
Suppress the index $i$ and consider $\kappa = \kappa_i$ and $\lambda_j = \lambda_{ij}$ as functions
of $(s_j) = (s_{ij})$. Due to invariance of the first differential, upper dots on the left hand side
of \eqref{eqn:SDotLDot} can be viewed as variations with respect to the variables $(s_j)$, and we
find ourselves in the setting of the present subsection.

By assumption, we have $\dot \kappa = 0$ for the variation $\dot s$ at the point $s^0$. Together
with equation \eqref{eqn:DotKappa} (or \eqref{eqn:DotKInS}) this implies
$$
\sum_j s^0_j \dot\lambda_j = - \dot \kappa = 0.
$$
On the other hand,
$$
\sum_j s^0_j \dot\lambda_j = \sum_{j,k} s^0_j \frac{\partial \lambda_j}{\partial s_k} \dot s_k =
\sum_{j,k} s^0_j \frac{\partial^2 K}{\partial s_j \partial s_k} \dot s_k = D_{s^0}^2 K (s^0, \dot
s).
$$
Thus we have
$$
D_{s^0}^2 K (s^0, \dot s) = 0,
$$
which means that the vectors $s^0$ and $\dot s$ are mutually orthogonal with respect to the
symmetric bilinear form $D_{s^0}^2 K$. By Lemma \ref{lem:SignD2K}, $s^0$ is a positive vector for
$D_{s^0}^2 K$, and $D_{s^0}^2 K$ is negative semidefinite on the orthogonal complement to $s^0$. It
follows that
$$
D_{s^0}^2 K (\dot s, \dot s) \le 0.
$$
Since, by \eqref{eqn:DDotKInS},
$$
\sum_j \dot s_j \dot \lambda_j = \ddot K = D_{s^0}^2 K (\dot s, \dot s),
$$
the inequality \eqref{eqn:SDotLDot} in Lemma \ref{lem:KDDotNegative} follows. Equality $\sum_j \dot
s_j \dot \lambda_j = 0$ means that $\dot s \in \ker D_{s^0}^2 K$. This implies $\dot \lambda_j = 0$
for all $j$, as established in the first lines of the proof of Lemma \ref{lem:SignD2K}.

The lemma is proved.
\end{proof}

\subsection{The case of a non-simplicial polyhedron $P$}
\label{subsec:NonSimplicial}
First of all, if $P$ has some non-triangular faces, our definition of an isometric infinitesimal
deformation (the second half of Definition \ref{dfn:InfDeform}) has to be modified. If we only
require stability of lengths of edges, then, say, the cube would be considered as infinitesimally
flexible. There are two possibilities. Either one requires that each face of $P$ is moved by $(q_i)$
as a rigid plate. Or, one subdivides each face into triangles by non-crossing diagonals and requires
stability of lengths also for the diagonals. The latter class of infinitesimal deformations is a
priori larger than the former, and can be shown to be independent of the choice of subdividing
diagonals.

We choose the second possibility. Thus, $P = P(r^0)$ can be viewed as a simplicial polyhedron with
some dihedral angles equal to $\pi$. As we vary $r$ in a neighborhood $U$ of $r^0$, these dihedral
angles can become less than $\pi$. Undaunted by this, we carry out our arguments for non-convex
warped polyhedra as well. As we come to warped spherical polygons, we consider also non-strictly
convex and ``slightly non-convex'' ones. Everything goes through, including Lemma \ref{lem:SignD2K}.
For the positivity of the form $D_{s^0}^2 K$ on the vector $s^0$ we need only positivity of the area
of the corresponding Euclidean polar dual. This is fulfilled, as the polygon $S$ is a convex (albeit
non-strictly) spherical polygon.

Thus a simple modification of the definitions allows to extend the argument to non-simplicial
polyhedra.

\section{Duality}
\subsection{Duality between second derivatives of $\Vol$ and $\HE$}
\label{subsec:VolHEDual}
There is an apparent duality between constructions and arguments in Sections \ref{sec:GaussRig} and
\ref{sec:MetrRig}. This is the same kind of duality as between Cauchy's proof of Theorem
\ref{thm:PolyhRig} in \cite{Cau13} and Alexandrov's proof of Theorem \ref{thm:MinkRigDiscr} in \cite[Chapter XI]{Ale05}. Cauchy studies a deformation of the spherical link $S_i$ of a vertex $p_i$, marks with $+$ or $-$ the $j$-th vertex of $S_i$ if the angle $(\lambda_{ij})$ increases or decreases during the deformation and shows that either at least four sign changes occur as one goes along the boundary of $S_i$, or no vertex of $S_i$ is marked, i.~e. the dihedral angles at $p_i$ are stable. Then a combinatorial argument involving the Euler characteristic of a polyhedron shows that there are no markings at all. Alexandrov considers variations of edge lengths along the boundary of a face $F_i$, marks lengthening and shortening edges, and follows Cauchy's argument.

In this paper, new similarities show up. Deformations of $Q$ in Section \ref{sec:GaussRig} are
governed by support parameters $(h_i)$, while deformations of $P$ in Section \ref{sec:MetrRig} are
governed by radii $(r_i)$. There are functions $\Vol(h)$ and $\HE(r)$ with similar variational
properties, and infinitesimal rigidity is proved by showing that
$$
\dim \ker (D^2 \Vol) = 3 \, \mbox{ and } \dim \ker (D^2 \HE) = 3.
$$
The following lemma strengthens this analogy by making it quantitative rather than qualitative.

\begin{lem}
\label{lem:DualEucl}
Let $P \subset \R^3$ be a convex polyhedron with vertices $p_1, \ldots, p_n$ and such that $0$ lies
in the interior of $P$. Let $Q = P^*$ be the polar dual of~$P$, that is
$$
Q = \{x \in \R^3\ |\ \langle x, \nu_i \rangle \le h_i^0, i=1,\ldots,n\},
$$
where
$$
\nu_i = \frac{p_i}{\|p_i\|} = \frac{p_i}{r_i^0}, \quad h_i^0 = \frac{1}{r_i^0}.
$$
Then we have
\begin{equation}
\label{eqn:Mystery}
D_{h^0}^2 \Vol = D_{r^0}^2 \HE.
\end{equation}
\end{lem}
\begin{proof}
A direct computation \cite[Section 3]{BI08} shows that
$$
\frac{\partial^2 \HE}{\partial r_i \partial r_j} = \frac{\cot \beta_{ij} + \cot
\gamma_{ij}}{\ell_{ij} \sin \rho_{ij} \sin \rho_{ji}} = \frac{\partial^2 \Vol}{\partial h_i \partial
h_j},
$$
$$
\frac{\partial^2 \HE}{\partial r_i^2} = - \sum_j \cos \phi_{ij}\frac{\cot \beta_{ij} + \cot
\gamma_{ij}}{\ell_{ij} \sin \rho_{ij} \sin \rho_{ji}} = \frac{\partial^2 \Vol}{\partial h_i^2}.
$$
Here angles $\beta_{ij}$ and $\gamma_{ij}$ were defined in the proof of Lemma \ref{lem:KFirstVar},
$\phi_{ij}$ is the angle at $0$ in the triangle $0p_ip_j$, and $\ell_{ij}$ and $\rho_{ij}$ are as on
Figure \ref{fig:SpherLink}. The lemma is proved.
\end{proof}

An immediate consequence of Lemma \ref{lem:DualEucl} is that each of Theorems \ref{thm:MinkRigDiscr}
and \ref{thm:PolyhRig} implies the other one (e.~g. one can replace Subsections \ref{subsec:DHE} to
\ref{subsec:SecVarK} by Lemma \ref{lem:DualEucl}). On the other hand, a proof of identity
\eqref{eqn:Mystery} by a direct computation is more a question than an answer. There is no obvious
reason why \eqref{eqn:Mystery} should hold, for example the functions $\HE(r^0 + x)$ and $\Vol(h^0 + x)$ are by no means equal.

The next subsection partially ``demystifies '' identity \eqref{eqn:Mystery} by pointing out a close
relationship between $\HE$ and $\Vol$ in the context of polar duality in the $3$--sphere.

\subsection{Duality in spherical geometry}
\label{subsec:SpherDual}
Let $P \subset \Sph^3$ be a convex spherical polyhedron (i.~e. the intersection of a finite number
of hemispheres that contains no pair of antipodal points). Define the polar dual of $P$ as
$$
P^* = \{x \in \Sph^3\ |\ \langle x, y \rangle \le 0 \mbox{ for all } y \in P\}.
$$
(We put $\Sph^3 = \{x \in \R^4\ |\ \|x\| = 1\}$, and let $\langle \cdot, \cdot \rangle$ be the
scalar product in $\R^4$.) Then $P^*$ is also a finite intersection of hemispheres (take one
hemisphere for every vertex of $P$). We assume that $\dim P = 3$; then $P^*$  contains no pair of
antipodal points. Let $p_ip_j$ be an edge of $P$. In $P^*$ there is a dual edge lying in the
intersection of hyperspheres polar to $p_i$ and $p_j$. Denote by $\ell_{ij}$ the length of the edge
$p_ip_j$, and by $\ell^*_{ij}$ the length of the dual edge.

\begin{thm}[P.~McMullen, \cite{McM75}]
\label{thm:WeylTube}
For every $3$--dimensional spherical polyhedron $P$ and its polar dual $P^*$ the following identity
holds:
\begin{equation}
\label{eqn:WeylTube}
\Vol(P) + \frac{1}{2} \sum_{\{i,j\}} \ell_{ij} \ell^*_{ij} + \Vol(P^*) = \pi^2.
\end{equation}
\end{thm}

The proof is based on Lemma \ref{lem:McM} below.

We have
$$
\ell_{ij} = \lambda^*_{ij}, \quad \ell^*_{ij} = \lambda_{ij},
$$
where $\lambda_{ij}$ is the exterior dihedral angle of $P$ at the edge $p_ip_j$, and
$\lambda^*_{ij}$ is the exterior dihedral angle of $P^*$ at the dual edge. The first equation holds
because the distance between $p_i$ and $p_j$ is equal to the exterior dihedral angle between the
hyperspheres polar to $p_i$ and $p_j$. The second holds since $(P^*)^* = P$ and $p_ip_j$ is dual of
its dual. It follows that the summand in the middle of \eqref{eqn:WeylTube} is the discrete total mean curvature of $\partial P$ and at the same time the discrete total mean curvature of $\partial P^*$.

More generally, let $P$ be a convex polyhedron in the $d$--sphere $\Sph^d$ for arbitrary $d \ge 1$, and let $F$ be a face of $P$. We define the dual face of $P^*$ as
$$
F^\perp = \{x \in P^*\ |\ \langle x, y \rangle = 0 \mbox{ for all } y \in F\}.
$$
We have $\dim F^\perp + \dim F = d - 1$. Every polyhedron is considered to be its own face and also to have $\emptyset$ as a $(-1)$--dimensional face. Then, clearly
$$
\emptyset^\perp = P^*, \quad P^\perp = \emptyset.
$$
Define the norm of a $k$-dimensional face as its $k$--volume divided by the $k$--volume of the
$k$-dimensional sphere:
$$
\|F\| = \frac{\Vol_k (F)}{\Vol_k (\Sph^k)}.
$$
As $\Vol_0(\Sph^0) = 2$, the norm of a vertex equals $\frac{1}{2}$; we also put $\|\emptyset\| = 1$.

\begin{figure}[ht]
\begin{center}
\input{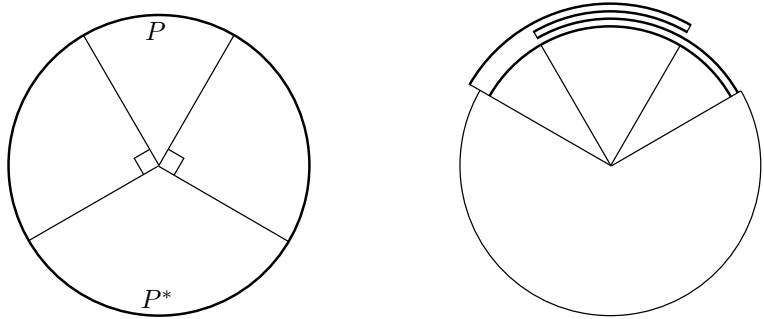}
\end{center}
\caption{To the proof of Lemma \ref{lem:McM}.}
\label{fig:Coverings}
\end{figure}

\begin{lem}
\label{lem:McM}
For every convex $d$--dimensional polyhedron $P \subset \Sph^d$ the following identities hold:
\begin{equation}
\label{eqn:Steiner}
\sum_F \|F\| \cdot \|F^\perp\| = 1,
\end{equation}
\begin{equation}
\label{eqn:SteinerNeg}
\sum_F (-1)^{\dim F} \|F\| \cdot \|F^\perp\| = 0.
\end{equation}
Here the summation extends over all faces of $P$, including $\emptyset$ and $P$ itself.
\end{lem}
\begin{proof}
Denote by $F \times F^\perp$ the convex hull of $F \cup F^\perp$ in $\Sph^d$. For every $F$, this is a
$d$--dimensional spherical polyhedron, and we have
$$
\|F \times F^\perp\| = \|F\| \cdot \|F^\perp\|,
$$
since $\span(F)$ and $\span(F^\perp)$ are mutually orthogonal subspaces of $\R^{d+1}$. It is not hard to
see that the sphere decomposes as union of polyhedra
$$
\Sph^d = \bigcup_F (F \times F^\perp),
$$
with disjoint interiors. Thus we have
$$
\Vol_d(\Sph^d) = \sum_F \Vol_d(F \times F^\perp)
$$
which implies equation \eqref{eqn:Steiner}.

Equation \eqref{eqn:SteinerNeg} is proved in a similar way, by replacing $P^*$ with $-P^*$. The
family of polyhedra $F \times -F^\perp$, each counted with multiplicity $(-1)^{\dim F}$, forms a
covering with total multiplicity $0$ over every point of $\Sph^d$. Therefore
$$
\sum_F (-1)^{\dim F} \Vol_d(F \times F^\perp) = 0,
$$
and equation \eqref{eqn:SteinerNeg} follows. Figure \ref{fig:Coverings} illustrates the case $d=1$.
\end{proof}

\begin{proof}[Proof of Theorem \ref{thm:WeylTube}]
Since $\Vol_3(\Sph^3) = 2\pi^2$ and $\Vol_1(\Sph^1) = 2\pi$, for $d=3$ formulas \eqref{eqn:Steiner} and \eqref{eqn:SteinerNeg} yield
\begin{gather*}
\Vol(P) + \frac{\pi}4 \Area(\partial P) + \frac12 \sum_{\{i,j\}} \ell_{ij} \ell^*_{ij} + \frac{\pi}4 \Area(\partial P^*) + \Vol(P^*) = 2\pi^2,\\
\Vol(P) - \frac{\pi}4 \Area(\partial P) + \frac12 \sum_{\{i,j\}} \ell_{ij} \ell^*_{ij} - \frac{\pi}4 \Area(\partial P^*) + \Vol(P^*) = 0.
\end{gather*}
By adding and subtracting these two formulas, we obtain
\begin{gather}
\Vol(P) + \frac{1}{2} \sum_{\{i,j\}} \ell_{ij} \ell^*_{ij} + \Vol(P^*) = \pi^2, \notag\\
\Area(\partial P) + \Area(\partial P^*) = 4\pi \label{eqn:GB}
\end{gather}
Theorem \ref{thm:WeylTube} is proved.
\end{proof}

\begin{rem}
Equation \eqref{eqn:GB} is the discrete Gauss-Bonnet theorem for surfaces in $\Sph^3$. Equation \eqref{eqn:Steiner} is the Steiner formula for the volume of the $\frac{\pi}{2}$--neighborhood of~$P$. Equation \eqref{eqn:WeylTube} is an analog of the Weyl tube formula.
\end{rem}

Lemma \ref{lem:McM} is proved by Peter McMullen in \cite{McM75}. Section 4 of \cite{McM75} contains
also references to smooth analogs of formulas \eqref{eqn:WeylTube} and \eqref{eqn:GB}. The main reference is Herglotz \cite{Her43b} (see also \cite[Subsection 6.4]{Izm11b}) who mentions that Fenchel has proved formulas \ref{eqn:WeylTube} and \eqref{eqn:GB} by induction on dimension. There seems to be no written account of Fenchel's proof.

\begin{rem}
\label{rem:SphInt}
Points of $P^*$ are poles of great spheres disjoint from the interior of $P$. This gives Theorem
\ref{thm:WeylTube} an integral-geometric interpretation: a random great sphere intersects $P$ with
the probability
$$
\frac{1}{\pi^2} \left( \Vol(P) + \frac12 \sum_{\{i,j\}} \ell_{ij} \lambda_{ij} \right).
$$
See also \cite[Chapter 17, \S 5, Note 1]{San04}.
\end{rem}

\begin{rem}
Milnor on the last page of \cite{Mil94} suggests another way of proving Equation \eqref{eqn:WeylTube}, by deforming $P$ to a point and integrating the Schl\"afli formula. It would be interesting to see whether the Schl\"afli formula can be derived from \eqref{eqn:WeylTube}.
\end{rem}

\subsection{Gauss and metric rigidity for convex spherical polyhedra}
\label{subsec:SpherDualMetr}
The following are the spherical analogs of Theorems \ref{thm:MinkRigDiscr} and \ref{thm:PolyhRig}.
\begin{thm}
\label{thm:MinkRigSph}
Let $Q \subset \Sph^3$ be a simple convex polyhedron, and let $(Q_t)$ be a deformation of
$Q$ that keeps all dihedral angles constant in the first order. Then $Q_t$ is congruent to $Q$ in
the first order.
\end{thm}
This theorem extends to non-simple convex spherical polyhedra, if one requires stability of dihedral
angles between the planes of any two faces having a vertex in common. A less restrictive
generalization allows change of combinatorics and appearance of negative edge lengths, see
Subsection~\ref{subsec:NonSimple}.

To substantiate the analogy between Theorems \ref{thm:MinkRigDiscr} and \ref{thm:MinkRigSph}, note
that directions of face normals of a Euclidean polyhedron determine its dihedral angles, and that
dihedral angles of a simple polyhedron determine angles of its faces, and thus, for a spherical
polyhedron, face areas.

\begin{thm}
\label{thm:PolyhRigSph}
Let $P \subset \Sph^3$ be a simplicial convex spherical polyhedron, and let $(P_t)$ be a deformation
of $P$ that keeps all edge lengths constant in the first order. Then $P_t$ is congruent to $P$ in
the first order.
\end{thm}
A generalization can be stated for non-simplicial polyhedra by subdividing non-triangular faces, see
Subsection \ref{subsec:NonSimplicial}.

Note that Theorems \ref{thm:MinkRigSph} and \ref{thm:PolyhRigSph} are equivalent: just put $Q_t =
P_t^*$.

One can prove Theorem \ref{thm:PolyhRigSph} following the arguments used in Section \ref{sec:MetrRig} to prove its Euclidean analog. Warped spherical polyhedra are defined in an obvious way. In place of the Hilbert-Einstein functional consider the functional
$$
S(P) := 2 \Vol(P) + \sum_i r_i \kappa_i + \sum_{\{i,j\}} \ell_{ij} \lambda_{ij},
$$
where $P = P(r)$ is a warped polyhedron with radii $(r_i)$. Schl\"afli's formula implies
$$
\frac{\partial S}{\partial r_i} = \kappa_i.
$$
Functional $S$ can be described as a discrete gravity action with non-zero cosmological constant, cf. \cite[Subsection 6.5]{Izm11b}.

When we want to copy the approach of Section \ref{sec:GaussRig} to prove Theorem
\ref{thm:MinkRigSph}, it undergoes more substantial changes and becomes more close to the approach to the dual theorem. Let $q_0$ be an interior point of $Q$. By dropping
perpendiculars $q_0q_i$ to the faces and then perpendiculars $q_i q_{ij}$ to the edges of $Q$, we
cut $Q$ into ``bricks''. A brick is a polyhedron combinatorially equivalent to the cube and with six
right angles at the edges not incident to the vertices $q_0$ and $q_{ijk}$.
A brick is uniquely determined by lengths $h_i, h_j, h_k$ of edges adjacent to $q_0$ and by dihedral
angles at the edges adjacent to $q_{ijk}$ (just choose $q_0$ as a unique point at distances $h_i,
h_j, h_k$ from the planes of a given trihedral angle). We vary ``support parameters'' $(h_i)$ of a
brick decomposition of $Q$ while leaving dihedral angles of $Q$ constant. As a result, singularities
around the edges $q_0q_i$ appear. Denote by $\kappa_i$ the curvatures at $q_0q_i$ and put
$$
S^*(Q) := 2 \Vol(Q) + \sum_i h_i \kappa_i.
$$
Then again Schl\"afli's formula implies
$$
\frac{\partial S^*}{\partial h_i} = \kappa_i.
$$
As deformations of $Q$ that keep dihedral angles constant can be identified with vectors $\dot h$
such that $\dot \kappa = 0$, Theorem \ref{thm:MinkRigSph} reduces to a statement about the kernel of $D^2 S^*$.

\begin{lem}
\label{lem:HESphDual}
Let $P$ be a warped spherical polyhedron, and $P^*$ be its polar dual. Then we have
$$
S(P) + S^*(P^*) = 2\pi^2.
$$
\end{lem}
\begin{proof}
Let $\Sph^3_\kappa$ be a spherical cone manifold obtained by taking a warped product of $[0,\pi]$ with the spherical link $S_0$ of the warped polyhedron $P$. (The manifold $\Sph^3_\kappa$ is
a $3$--sphere with $n$ singular meridians of curvatures $(\kappa_i)$; meridians are arranged as the cone singularities in the link of the point $p_0 \in P$.) Both warped polyhedra $P$ and $P^*$ can be embedded in $\Sph^3_\kappa$.

Similarly to Lemma \ref{lem:McM} and Theorem \ref{thm:WeylTube}, we can show that
$$
\Vol(P) + \frac12 \sum_{\{i,j\}} \ell_{ij} \ell^*_{ij} + \Vol(P^*) = \frac\pi{2} (2\pi - \sum_i
\kappa_i),
$$
the right hand side being half the volume of $\Sph^3_\kappa$. On the other hand, radii of $P$ and
heights of $P^*$ are related by $h_i = \pi - r_i$. Therefore
\begin{multline*}
S(P) + S^*(P^*) = 2 \Vol(P) + \sum_{\{i,j\}} \ell_{ij} \ell^*_{ij} + 2 \Vol(P^*) +
\sum_i (r_i + h_i)\kappa_i\\
= \pi (2\pi - \sum_i \kappa_i) + \pi \sum_i \kappa_i = 2\pi^2,
\end{multline*}
and the lemma is proved.
\end{proof}

Thus in the spherical case the duality between Gauss and metric rigidity of convex polyhedra is perfect, and their proofs that use the second derivatives of $S^*$ and $S$ are simply the same.

\subsection{Duality in hyperbolic-de Sitter geometry}
\label{subsec:HypDS}
Consider the hyperboloid model of the hyperbolic space
$$
\H^3 = \{x \in \R^4\ |\ \|x\|_{3,1} = -1, x_0 > 0\},
$$
where $\|x\|_{3,1} = -x_0^2 + x_1^2 + x_2^2 + x_3^2$ is the Minkowski scalar product. The orthogonal complement $x^\perp \subset \R^4$ does not intersect $\H^3$, but it intersects the one-sheeted hyperboloid
$$
\dS^3 = \{x \in \R^4\ |\ \|x\|_{3,1} = 1\}
$$
which is a model of the \emph{de Sitter space}. For every convex polyhedron $P \subset \H^3$, define its polar dual as
$$
P^* = \{x \in \dS^3 \cup \H^3_- \ | \langle x, y \rangle_{3,1} \le 0\},
$$
where $\H^3_-$ is the antipodal copy of $\H^3$. In particular, $\partial P^*$ is a convex polyhedral surface in $\dS^3$. Hyperbolic-de Sitter duality was used by Rivin \cite{Riv86} to study geometry of hyperbolic polyhedra.

The Minkowski scalar product induces a semi-Riemannian metric on $\dS^3$, and the boundary of $P^*$
is space-like with respect to this metric. Between lengths and angles in polar
objects in $\H^3$ and $\dS^3$ there is a correspondence similar to that in $\Sph^3$, see Subsection
\ref{subsec:SpherDual}. As a consequence, Gauss (respectively, metric) rigidity of a convex
polyhedron in $\H^3$ is equivalent to metric (respectively, Gauss) rigidity of its polar dual in
$\dS^3$.

A variational approach to the metric rigidity of convex hyperbolic polyhedra is based on the
discrete gravity action
$$
S(P) := - 2 \Vol(P) + \sum_i r_i \kappa_i + \sum_{\{i,j\}} \ell_{ij} \lambda_{ij}
$$
where $P$ is a warped hyperbolic polyhedron. Similarly to the previous Subsection, we have
$$
S(P) + S^*(P^*) = 0,
$$
where $S^*$ is the gravity action without the boundary term. This is implied by the following analog of Theorem \ref{thm:WeylTube}.

\begin{thm}
\label{thm:WeylTubeHyp}
For every convex polyhedron $P \subset \H^3$ and its polar dual $P^* \subset \dS^3 \cup \H^3_-$ the following
identity holds:
\begin{equation}
\label{eqn:WeylTubeHyp}
\Vol(P) - \frac{1}{2} \sum_{\{i,j\}} \ell_{ij} \ell^*_{ij} + \Vol(P^*) = 0.
\end{equation}
\end{thm}

Note that $P^*$ is the union of a whole hyperbolic space $\H^3_-$ and of an infinite end of $\dS^3$. However, there is a consistent way to define a finite measure $\Vol(P^*)$, \cite{BIH92, CK06}.

Similarly to Milnor's approach in the spherical case, Theorem \ref{thm:WeylTubeHyp} is proved in \cite{SP00} by integrating the Schl\"afli formula.

\begin{rem}
Formula \eqref{eqn:WeylTubeHyp} has the following integral-geometric interpretation: the (motion-invariant and appropriately normalized) measure of the set of all planes that intersect a convex body $P \subset \H^3$ equals
$$
\frac{1}{2} \sum_{\{i,j\}} \ell_{ij} \lambda_{ij} - \Vol(P).
$$
In particular, this quantity is monotone under inclusion and always positive.

Integral non-euclidean geometry is dealt with in \cite[Chapter 17]{San04}. A proof of Theorems \ref{thm:WeylTube} and \ref{thm:WeylTubeHyp} indicated in \cite[Chapter 17, \S 5, Note 1]{San04} is by obtaining first a similar formula in the smooth case and then going to the limit $\epsilon \to 0$ for the boundary of $\epsilon$--neighborhood of $P$.
\end{rem}

\subsection{Shearing vs.~bending}
\label{subsec:ShearBend}
There is another kind of duality between Theorems \ref{thm:MinkRigDiscr} and \ref{thm:PolyhRig},
relating metric rigidity of $P$ with Gauss rigidity of~$P$, rather than with that of
$P^*$.

An infinitesimal isometric deformation of $P$ can be described by assigning to each face of $P$ an
infinitesimal isometry of $\R^3$, so that these isometries agree on the edges. (Faces move as rigid
plates joined by hinges.) An infinitesimal isometry of $\R^3$ is an infinitesimal screw motion, that is a vector field
$$
\xi(x) = \eta \times x + \tau, \quad \eta, \tau \in \R^3.
$$
Let $\xi_i = (\eta_i, \tau_i)$ be the infinitesimal isometry associated to a face $F_i$. Let $F_{ij} = F_i \cap F_j$ be an edge of $P$. The condition $\xi_i|_{F_{ij}} = \xi_j|_{F_{ij}}$ is equivalent to
$$
(\eta_i - \eta_j) \times x = \tau_j - \tau_i \quad \mbox{ for all } x \in F_{ij},
$$
which implies that the vector $\eta_i - \eta_j$ is parallel to $F_{ij}$:
\begin{equation}
\label{eqn:DPar}
\eta_i - \eta_j \parallel F_{ij} \quad \mbox{ for all edges } F_{ij}.
\end{equation}

Now, instead of rotating each face $F_i$ according to the vector $\eta_i$, let us translate the plane of $F_i$ by $\eta_i$ (consider a plane as a set of points rather than as a geometric figure, so that the horizontal component of $\eta_i$ cannot be neglected). By \eqref{eqn:DPar}, the two translations of $\span(F_{ij})$ differ by a vector parallel to $F_{ij}$. This implies that the boundary of $F_i$ is translated together with $F_i$, except that some new edges may appear. Namely, if $F_i$ and $F_k$ had only a vertex in common, the translations can split this vertex and create an edge between $F_i$ and $F_k$. If there is no vertex splitting at all, then the polyhedron $P$ is translated as a rigid body.

Even if vertex splitting happens, the area of each face changes by $o(t)$, if we do translations by
$t\eta_i$. It follows that putting $\dot h_i = \langle \eta_i, \nu_i \rangle$ we obtain $\dot A_i = 0$ for all $i$, in the notation of Subsection \ref{subsec:Theorem}.

Thus to every isometric infinitesimal deformation of $P$ there corresponds an infinitesimal Gauss
image preserving deformation of $P$. The correspondence can be inverted, by solving the above
equations for $(t_i)$.

This observation was made by Weyl in \cite{Weyl17} who was inspired by a work of Blaschke
\cite{Bla12} where a similar correspondence for smooth surfaces was indicated. See also
\cite[Chapter XI, \S 3]{Ale05}. A similar correspondence is known for ideal hyperbolic polyhedra, where it is sometimes described as the duality between shearing and bending. Indeed, with $(d_i)$ viewed as rotations, the difference $d_i - d_j$ represents the change of the dihedral angle between $F_i$ and $F_j$, whereas if $(d_i)$ are viewed as translations, then $d_i - d_j$ is the shift of $F_i$ along $F_{ij}$ relative to $F_j$.

\subsection{Statics and polarity}
\label{subsec:StaPol}
Infinitesimal rigidity of a polyhedron is equivalent to static rigidity of its 1--skeleton (with
diagonals of non-triangular faces added), see \cite{RW81,Izm09}. Basically, statics deals with duals
of vector spaces of isometric, respectively trivial, infinitesimal deformations.

In \cite[Theorems 3.1 and 3.2]{Whi87}, Whiteley establishes a correspondence between statics of a
$3$--dimensional polyhedron and of its polar dual. The polarity can be taken with respect to any
quadric. This reflects the projective nature of static (and hence infinitesimal) rigidity, see next
subsection.

\subsection{Projective invariance and infinitesimal Pogorelov maps}
Infinitesimal rigidity of a bar-and-joint framework in $\R^d$ is invariant under projective
transformations. This fact follows from a projective formulation of statics, see \cite{CW82,Izm09}.

Projective invariance provides a link between infinitesimal rigidity in $\R^d$, $\Sph^d$, and
$\H^d$: an infinitesimally rigid framework in $\R^d$ remains infinitesimally rigid when viewed as a
framework in a projective model of $\Sph^d$ or $\H^d$. The arising correspondences between
infinitesimal isometric deformations are called infinitesimal Pogorelov maps, see
\cite{Pog73,Izm09}.

\section{Miscellaneous remarks}
\subsection{Regge action and Steiner formula}
The sum $\sum_i r_i \kappa_i$ appearing in the definition of the discrete Hilbert-Einstein
functional for warped polyhedra makes sense for an arbitrary compact closed manifold glued from
Euclidean simplices. It was introduced in \cite{Reg61} and is sometimes called the Regge action. It
was shown in \cite{CMS84} that $\sum_i r_i \kappa_i$ converges to a constant times the total scalar
curvature if a sequence of piecewise Euclidean manifolds converges in some good sense to a smooth
manifold.

The sum $\frac12 \sum_{\{i,j\}} \ell_{ij} \lambda_{ij}$ is a discrete analog of the total mean
curvature. It appears as a coefficient at $\epsilon^2$ in the formula for the volume of an
$\epsilon$--neighborhood of a convex polyhedron (as the total mean curvature does for convex bodies
with smooth boundary). This was noticed by Jakob Steiner \cite{Ste40}. The convergence of
$\sum_{\{i,j\}} \ell_{ij} \lambda_{ij}$ to a constant times the total mean curvature follows from
the continuity of mixed volumes with respect to the Hausdorff distance, see \cite{Scn93}.

\subsection{Mixed volumes and the Alexandrov-Fenchel inequality}
\label{subsec:MixVol}
Our constructions in Section \ref{sec:GaussRig} are closely related to the theory of mixed volumes,
see e.~g. \cite{Scn93}. In particular, the proof of Theorem \ref{thm:MinkRigDiscr} is copied from (a
part of) the proof of Alexandrov-Fenchel inequalities, \cite{Ale37II,Ale96,Scn93}. The relation
between mixed volumes and derivatives of the volume is apparent from the formula
$$
\Vol(h+tk) = \Vol(h) + 3t \Vol(h,h,k) + 3t^2 \Vol(h,k,k) + t^3 \Vol(k)
$$
that may serve as the definition of mixed volumes $\Vol(h,h,k)$ and $\Vol(h,k,k)$. (One has to
assume that polyhedra $Q(h)$ and $Q(k)$ are combinatorially isomorphic, in order that addition of
support parameters $h$ and $k$ correspond to Minkowski addition of polyhedra.) Right hand sides of
the formulas in Remark \ref{rem:SymmProp} give alternative expressions for mixed volumes.

To prove Theorem \ref{thm:MinkRigDiscr}, it is sufficient to show that the symmetric bilinear form
$D^2\Vol$ has corank $3$. The proof of Alexandrov-Fenchel inequalities goes further. It establishes
that $D^2\Vol$ has signature $(+_1, 0_3, -_{n-4})$, and in general shows by induction that the
signature is $(+_1, 0_d, -_{n-d-1})$ for $d$--dimensional polyhedra with $n$ facets.

\subsection{Existence theorems}
\label{subsec:ExistThms}
Infinitesimal rigidity can sometimes be used to prove existence theorems. Assume that an object is
infinitesimally rigid with respect to some parameters. If the space of objects and the space of
parameters have equal dimension, then the inverse function theorem provides us with local existence
and uniqueness of an object with given parameters. In other words, the map $\{\mbox{object}\}
\mapsto \{\mbox{parameters}\}$ is a local homeomorphism. A local homeomorphism with good topological properties is a global homeomorphism, which implies that any given set of parameters defines a unique object.

This method is used in \cite{Ale05} to prove the following two theorems.
\begin{thm}[Minkowski theorem]
Let $\nu_1, \nu_2, \ldots, \nu_n \in \R^3$ be unit vectors spanning $\R^3$, and let $C_1, C_2,
\ldots, C_n$ be positive numbers such that
$$
\sum_i C_i \nu_i = 0.
$$
Then there is a unique convex polyhedron in $\R^3$ with outer face normals $(\nu_i)$ and respective
face areas $(C_i)$.
\end{thm}

\begin{thm}[Alexandrov theorem]
Let $g$ be a Euclidean cone-metric on~$\Sph^2$ (every such metric can be obtained by gluing a set of Euclidean triangles). Assume that the angles at all cone points are less than $2\pi$. Then
there is a unique convex polyhedron in $\R^3$ such that $g$ is an intrinsic metric of its boundary.
\end{thm}
Note that in both theorems the combinatorial structure (which pairs of faces are adjacent, which
pairs of cone points are joined by edges) is not given in advance and is practically impossible to
determine without finding the polyhedron in question.

Properties of the functionals $\Vol$ and $\HE$ suggest a variational approach to both theorems.
Indeed,
$$
\frac{\partial}{\partial h_i}\left(\Vol(h) - \sum_i h_i C_i\right) = A_i(h) - C_i
$$
implies that the polyhedron in the Minkowski theorem is a critical point of the function $\Vol(h) -
\sum_i h_i C_i$ on the set of polyhedra with outer face normals $(\nu_i)$. Similarly,
$$
\frac{\partial \HE}{\partial r_i} = \kappa_i
$$
implies that the polyhedron in the Alexandrov theorem is a critical point of $\HE$ on the set of
warped polyhedra. As stated in Subsection \ref{subsec:MixVol}, the second differential $D^2\Vol$ has
only one positive eigenvalue. By a happy coincidence, $\Vol$ is concave on the hyperplane $\{\sum_i
h_i C_i = 1\}$. At the maximum point (whose existence can be shown by a simple trick) face areas
$A_i$ are proportional to $C_i$, so that a scaling yields the desired polyhedron. This proof is due
to Minkowski himself, \cite{Min03}, see also \cite[Chapter VII, \S 2]{Ale05}.

The situation is more complicated with $\HE$ and the Alexandrov theorem. The discrete
Hilbert-Einstein functional on warped polyhedra is neither concave nor convex, moreover the
signature of its second differential is non-constant. Nevertheless, similarly to Subsection
\ref{subsec:SecVarK} one can show that $\ker D^2 \HE = \{0\}$, provided that $\kappa_i > 0$ and the
areas of Euclidean polar duals of the vertex links are positive. This allows to apply the inverse
function theorem to the map $r \mapsto \kappa$ and to find the desired polyhedron by constructing a
family of warped polyhedra with $\kappa \to 0$. This proof is given in \cite{BI08}.

Yuri Volkov, a student of Alexandrov, studied the discrete Hilbert-Einstein functional with the aim
to find a variational proof of Alexandrov's theorem similar to that of Minkowski theorem. His proof in \cite{Vol55, VP71}) doesn't use the function $\HE$, but is similar in the spirit. The polyhedron is found by minimizing the sum of radii $(r_i)$ over all warped polyhedra (this time with apex at a boundary vertex) with negative $(\kappa_i)$. See also \cite{Vol68} which is reprinted in the Appendix to \cite{Ale05}, where Volkov studies $\HE$ in order to give some a priori bounds for the embedding problem.

In some theorems of Alexandrov or Minkowski type, the functional $\HE$, respectively $\Vol$, happens to be concave. This is the case with the Alexandrov convex cap theorem \cite{Izm08}, and with the Alexandrov and Minkowski-type theorems for convex hyperbolic cusps \cite{FI09,FI}. However, in all these cases proofs of existence of a maximum point are quite difficult.

\subsection{Related work}
Our proof of infinitesimal rigidity of convex polyhedra in Section \ref{sec:MetrRig} is related to
works of Pogorelov and Volkov \cite{Pog56,Vol56} and of Schlenker \cite{Scl07}. In \cite[Section
3]{Scl07} an alternative proof of Lemma \ref{lem:SignD2K} can be found.

Duality between metric rigidity of $P$ and Gauss rigidity of $P^*$ involving theory of mixed volumes appears in Paul Filliman's work \cite{Fil92}. Lee generalized this duality \cite{Lee96}
to higher dimensions, relating it with Peter McMullen's theory of weihghts on polytopes
\cite{McM96}. See also the work of Tay, White, and Whiteley \cite{TWW95}.

The approach of Section \ref{sec:MetrRig} was extended by Schlenker to star-shaped polyhedra with
vertices in convex position, \cite{Scl09}. A further development of this technique is given in
\cite{IS10}. There the signature of $D^2\HE$ for an arbitrary triangulation of a convex polyhedron
is determined.

\subsection{Directions for future research}
\begin{prb}
Can equation \eqref{eqn:Mystery} be obtained as a limiting case of Theorem~\ref{thm:WeylTube}?
\end{prb}

\begin{prb}
Give a proof of Theorem \ref{thm:WeylTubeHyp} similar to McMullen's proof of Theorem \ref{thm:WeylTube}. Probably, all prerequisites are contained in \cite{CK06}. See also \cite{Her43b}.
\end{prb}

Lemma \ref{lem:DualEucl} relates metric rigidity of $P$ and Gauss rigidity of $P^*$, Subsection
\ref{subsec:StaPol} relates metric rigidity of $P$ and $P^*$, and Subsection \ref{subsec:ShearBend}
relates Gauss and metric rigidity of $P^*$. In the smooth case, this forms a part of Darboux's wreath, \cite[Subsection 5.3]{Izm11b}. Sauer \cite{Sau55} found a discrete analog of Darboux's wreath for polyhedral surfaces with quadrilateral faces and four-valent vertices.

\begin{prb}
Is there a Darboux's wreath for polyhedral surfaces of arbitrary combinatorics?
\end{prb}

Problems \ref{prb:CW}--\ref{prb:Che} are discussed in more detail in Subsections \ref{subsec:HypMan}--\ref{subsec:DiscrBochner}.

\begin{prb}
\label{prb:CW}
Reprove the Calabi-Weil rigidity of compact closed hyperbolic $3$--manifolds by showing
\begin{equation}
\label{eqn:DimKerHE}
\dim \ker \left( \frac{\partial \kappa_{ij}}{\partial r_{kl}} \right) = 3n,
\end{equation}
where $n$ is the number of vertices of a triangulation of a hyperbolic manifold.
\end{prb}

\begin{prb}
In a similar way, prove the infinitesimal rigidity of compact hyperbolic $3$--manifolds with convex polyhedral boundary.
\end{prb}

\begin{prb}
\label{prb:HK}
In a similar way, reprove the infinitesimal rigidity of compact hyperbolic cone-manifold with all cone angles less than $2\pi$.
\end{prb}

\begin{prb}
\label{prb:Che}
Reprove Cheeger's vanishing theorem with a discrete Bochner technique.
\end{prb}

The next problem generalizes problems \ref{prb:CW} and \ref{prb:HK} to higher dimensions.

\begin{prb}
Define discrete Einstein manifolds in dimensions greater than $3$ and prove their infinitesimal rigidity under suitable assumptions.
\end{prb}

To approach the above problem, one can start from the following consideration. Gauss infinitesimal rigidity holds for convex polyhedra in all dimensions. Convex polyhedra of higher dimensions are also infinitesimally rigid, in fact they are \emph{too rigid}. An alternative extension of the infinitesimal rigidity theorem to higher dimensions could be infinitesimal rigidity in the class of discrete Einstein metrics with a given restriction to the boundary. To define discrete Einstein metrics for a star-like triangulation of $P$, consider infinitesimal deformations $\dot h = u$ of $P^*$ that preserve volumes of facets in the first order. Putting $\dot r = u$ should define an infinitesimal Einstein deformation of $P$.

\begin{prb}
Define infinitesimal Einstein deformations of warped convex polyhedra in dimension $d > 3$ and prove Einstein infinitesimal rigidity of convex polyhedra.
\end{prb}

\appendix
\section{Some trigonometry}
The first two lemmas concern a spherical triangle with side lengths $a$, $b$, $c$, and values
$\alpha$, $\beta$, $\gamma$ of the
respective opposite angles.

\begin{lem}
\label{lem:DADAlpha}
The partial derivatives of the angles as functions of side lengths are:
\begin{eqnarray*}
\frac{\partial \alpha}{\partial a} & = & \frac1{\sin b \sin \gamma} = \frac1{\sin c \sin \beta},\\
\frac{\partial \alpha}{\partial b} & = & - \frac{\cot \gamma}{\sin b}.
\end{eqnarray*}
\end{lem}
\begin{proof}
Straightforward calculation using spherical cosine and sine laws.
\end{proof}

\begin{figure}[ht]
\begin{center}
\input{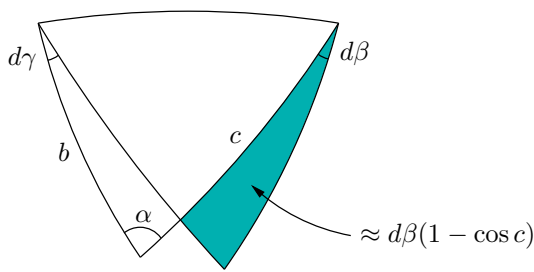}
\end{center}
\caption{Variation of the area of a spherical triangle when one of the side lengths is preserved.}
\label{fig:DotArea}
\end{figure}

\begin{lem}
\label{lem:DotABC}
Let the side lengths $b$ and $c$ change with velocities $\dot b$ and $\dot c$, while the
length $a$ remain constant. Then the corresponding variations of the angles satisfy the equation
\begin{equation}
\label{eqn:DotABC}
\dot\alpha + \dot\beta \cos c + \dot\gamma \cos b = 0.
\end{equation}
\end{lem}
\begin{proof}
This can be proved by tedious computations using Lemma \ref{lem:DADAlpha}. We give an alternative
geometric argument.

Denote by $\area$ the area of the triangle. As $\area = \alpha + \beta + \gamma - \pi$, equation
\eqref{eqn:DotABC} is equivalent to
$$
\dot \area = \dot \beta (1 - \cos c) + \dot \gamma (1 - \cos b).
$$
The latter is proved on Figure \ref{fig:DotArea} (recall that the area of a spherical cap of radius
$c$ equals $2\pi (1 - \cos c)$).
\end{proof}


\end{document}